\documentclass[twocolumn]{autart} 

\usepackage{graphicx}
\usepackage{amsmath}
\usepackage{amssymb}
\usepackage{mathtools} 
\usepackage[utf8]{inputenc}
\usepackage{letltxmacro}
\usepackage[table,xcdraw]{xcolor}
\usepackage{hhline} 

\newcommand{\nref}[1]{(\ref{#1})}
\newcommand{\bfs}[1]{\boldsymbol{#1}}
\newcommand{\col}[1]{\operatorname{col}\left({#1}\right)}
\newcommand{\diag}[1]{\operatorname{diag}\left({#1}\right)}
\newcommand{\kraj}{\hspace*{\fill} \qed}
\newcommand{\krajdokaz}{\hfill $\blacksquare$}

\newcommand{\z}[1]{\left( #1 \right)}

\newcommand{\m}[1]{\begin{bmatrix} #1 \end{bmatrix}}
\newcommand{\sm}[1]{\left[\begin{smallmatrix} #1 \end{smallmatrix}\right]}
\newcommand{\n}[1]{\left\| #1 \right\|}
\newcommand{\red}{\nonumber \\}
\newcommand{\vprod}[2]{\left \langle #1 \ \middle\vert\  #2 \right \rangle}
\newcommand{\oln}[1]{\overline{#1}}

\newcommand{\R}{\mathbb{R}}

\newcommand{\dom}{\mathrm{dom}}

\usepackage[utf8]{inputenc}
\usepackage[shortlabels]{enumitem} 
\usepackage{cuted}

\newcommand{\Diag}[1]{\operatorname{Diag}\left({#1}\right)}
\newcommand{\alphadonjii}{\underline{\alpha_{2,\rho}}\z{\n{x}_{\mathcal{M}_{\rho}}}}
\newcommand{\alphadonji}{\underline{\alpha_{1}}\z{\n{x_1}_\mathcal{A}}} 
\newcommand{\alphagornjii}{\overline{\alpha_{2,\rho}}\z{\n{x}_{\mathcal{M}_{\rho}}}}
\newcommand{\alphagornji}{\overline{\alpha_{1}}\z{\n{x_1}_\mathcal{A}}} 
\newcommand{\mrho}{{\mathcal{M}_{\rho}}}
\newcommand{\mA}{{\mathcal{M}_{\mathcal{A}}}}
\newcommand{\uajz}[1]{\underline{\alpha_{1}}\z{#1}} 
\newcommand{\uaj}{\underline{\alpha_{1}}} 
\newcommand{\oajz}[1]{\overline{\alpha_{1}}\z{#1}} 
\newcommand{\oaj}{\overline{\alpha_{1}}} 
\newcommand{\uadz}[1]{\underline{\alpha_{2,\rho}}\z{#1}} 
\newcommand{\uad}{\underline{\alpha_{2,\rho}}} 
\newcommand{\oadz}[1]{\overline{\alpha_{2,\rho}}\z{#1}} 
\newcommand{\oad}{\overline{\alpha_{2,\rho}}} 
\newcommand{\cA}{\mathcal{A}} 
\newcommand{\cX}{\mathcal{X}} 

\newtheorem{assum}{Assumption}
\newtheorem{theorem}{Theorem}
\newtheorem{lemma}{Lemma}
\newtheorem{remark}{Reamrk}

\newtheorem{corollary}{Corollary}
\newtheorem{defn}{Definition}
\newtheorem{example}{Example}

\allowdisplaybreaks 

\begin{document}

\begin{frontmatter}
\runtitle{Insert a suggested running title}  

\title{Stability of singularly perturbed hybrid systems with restricted systems evolving on boundary layer manifolds} 

\thanks[footnoteinfo]{This work was partially supported by the ERC under research project COSMOS (802348). E-mail addresses: \texttt{\{s.krilasevic-1, s.grammatico\}@tudelft.nl}.}

\author{Suad Krila\v sevi\' c} and
\author{\ Sergio Grammatico} 
          
\address{Delft Center for Systems and Control, TU Delft, The Netherlands}


\begin{abstract}                          

We present a singular perturbation theory applicable to systems with hybrid boundary layer systems and hybrid reduced systems {with} jumps from the boundary layer manifold. First, we prove practical attractivity of an adequate attractor set for small enough tuning parameters and sufficiently long time between almost all jumps. Second, under mild conditions on the jump mapping, we prove semi-global practical asymptotic stability of a restricted attractor set. Finally, for certain classes of dynamics, we prove semi-global practical asymptotic stability of the restricted attractor set for small enough tuning parameters and sufficiently long period between almost all jumps {of the slow states} only. 
\end{abstract}

\end{frontmatter}

\section{Introduction}

A realistic modeling of many control systems requires high-order nonlinear differential equations that might be difficult to fully analyze. To alleviate this problem, we often design control systems with various parameters that with proper tuning can effectively reduce the order of the model and thus simplify the  stability analysis. The main theoretical framework for such analysis is singular perturbation theory \cite{naidu1988singular}, \cite{khalil2002nonlinear}. The associated model reduction is accomplished by splitting the states into fast and slow states; for each constant value of the slow states, the fast states should converge to an equilibrium point defined by the slow states, and the union of these equilibrium points for all possible slow states defines the so-called boundary layer manifold. Then, the reduced system contains just the slow states and their dynamics assuming they are evolving along that manifold. \\ \\

Singular perturbation theory has been successfully applied to equilibrium seeking in optimization and game theory. One common method of applying zeroth-order algorithms to dynamical systems with cost measurements as output is through a time-time scale separation of the controller and the plant, as demonstrated in \cite{krstic2000stability} and \cite{poveda2017framework}. Time-scale separation can be useful for algorithms where consensus on specific states must be reached before initiating the equilibrium seeking process \cite{Carnevale2023nonconvex}, \cite{Ochoa2022momentum}, \cite{Wang2023distributed}, \cite{Sun2021Continous}. Furthermore, in some works \cite{poveda2017framework}, \cite{poveda2019robust}, \cite{poveda2020fixed}, via singular perturbation analysis the (pseudo)gradient estimate is filtered before being incorporated into the algorithm. Singular perturbation theory is also used to demonstrate algorithm convergence in problems with slowly varying parameters \cite{Galarza2022Sliding}.\\ \\

Several extensions of singular perturbation theory are known for hybrid systems. In \cite{sanfelice2011singular}, the authors examine a singularly perturbed system in which the boundary layer system is continuous, and the reduced system is hybrid, and both render the corresponding sets globally asymptotically stable. While the work in \cite{wang2012analysis} proposes averaging theory results, in can also be used to prove stability in singularly perturbed systems. Similarly to \cite{sanfelice2011singular}, the authors assume that the boundary layer system is continuous and that the averaged system, which plays the role of the reduced system, is hybrid. In \cite{wang2012averaging}, the same authors extend the results for the case when the boundary layer system itself is hybrid. In the aforementioned works, the reduced system is derived by assuming that the slow states ``flow" along the boundary layer manifold, while the slow states \emph{do not} jump from that manifold. Therefore, the reduced system jumps cannot use the properties of the boundary layer manifold to support stabilization; essentially only the continuous dynamics are used to prove stability, ``despite" the jumps.\\
In order for the discrete-time dynamics to support stabilization of singularly perturbed systems, we can design the dynamics so that we jump when we are in the proximity of the manifold. This scenario in principle is similar to that in \cite{teel2001solving}, \cite{khong2013unified} where the authors prove that there exists a sampling period such that a discrete-time optimization-based controller (the reduced system) can find a neighborhood of the optimum of a steady-state output map of a continuous system with an input (boundary layer system). In \cite{poveda2017robust}, the authors take a step further and design an event-triggered framework to accomplish the same task by measuring the changes in the output and in turn to determine when the system has approached the boundary layer manifold. Although these methods better incorporate discrete-time reduced system dynamics, the boundary layer system is still only continuous. In this paper, we instead deal with a hybrid boundary layer system and thus extend the current state of the art.\\
\emph{Contribution}: In view of the  above  literature, our theoretical contributions are summarized next:
\begin{itemize}
    \item We propose a singular perturbation theory for hybrid systems, where the reduced system takes into account the jumps from the boundary layer manifold, differently from \cite{sanfelice2011singular}, \cite{wang2012analysis} where jumps are assumed not to interfere with stability. Furthermore, we allow for the set of fast variables not to be bounded a priori, thus enabling the use of reference trajectories and counter variables in the boundary layer system.
    \item We prove semi-global practical asymptotic stability of the restricted attractor set, under certain mild assumptions on the jump mapping. This attractor set includes only the steady-state values of the fast states that correspond to the slow attractor states, rather than the complete range of possible fast variables.
    \item We show that, in a system resembling the one described in \cite{wang2012averaging}, where a distinction is made between jumps in the slow and fast states, the aforementioned results remain valid if there are sufficiently long intervals between nearly all jumps in the \emph{slow states}. 
\end{itemize}

Our theory enables the analysis of multiple timescale control systems where both the controller and the plant are hybrid. Furthermore, as the jumps occur at the boundary layer, it would be also possible to incorporate state/output feedback into the controller jump mappings.\\ \\

\emph{Notation}: {The set of real numbers and the set of nonnegative real numbers are denoted by $\mathbb{R}$ and  $\mathbb{R}_+$, respectively. Given a set $\mathcal{Z}$, $\mathcal{Z}^n$ denotes the Cartesian product of $n$ sets $\mathcal{Z}$.} For vectors $x, y \in \mathbb{R}^{n}$ and $\mathcal{A} \subset \R^n$, $\vprod{x}{y}$, $\|x \|$ and $\|x \|_\mathcal{A}$ denote the Euclidean inner product, norm, weighted norm and distance to set respectively. Given $N$ vectors $x_1, \dots, x_N$, possibly of different dimensions, $\col{x_1, \dots x_N} \coloneqq \left[ x_1^\top, \dots, x_N^\top \right]^\top $. Collective vectors are denoted in bold, i.e,  $\bfs{x} \coloneqq \col{x_1, \dots, x_N}$ as they collect vectors from multiple agents. We use $\mathbb{S}^{1}:=\left\{z \in \mathbb{R}^{2}: z_{1}^{2}+z_{2}^{2}=1\right\}$ to denote the unit circle in $\R^2$. $\operatorname{Id}$ is the identity operator; $I_n$ is the identity matrix of dimension $n$ and $ \bfs{0}_n$ is vector column of $n$ zeros; {their index is omitted where the dimensions can be deduced from context. The unit ball of appropriate dimensions depending on context is denoted with $\mathbb{B}$.} A continuous function $\gamma: \R_+ \rightarrow \R_+$ is of class $\mathcal{K}$ if it is zero at zero and strictly increasing. A continuous function $\alpha: \R_+ \rightarrow \R_+$ is of class $\mathcal{L}$ if is non-increasing and converges to zero as its arguments grows unbounded. A continuous function $\beta: \R_+ \times \R_+ \rightarrow \R_+$ is of class $\mathcal{KL}$ if it is of class $\mathcal{K}$ in the first argument and of class $\mathcal{L}$ in the second argument. UGAS refers to uniform global asymptotic stability, as defined in \cite[Def. 2.2, 2.3]{poveda2021robust}. We define semi-global practical asymptotic stability (SGPAS) similarly as in \cite{teel1999semi}.
\begin{defn}[SGPAS]
    The set $\mathcal{A}$ is SGPAS as $(\varepsilon_1, \dots, \varepsilon_k) \rightarrow 0$ for the parametrized hybrid system $\mathcal{H}_\varepsilon$, if for each given $\Delta > \delta > 0$,  there exists a parameter $\varepsilon_1^*$ such that for each $\varepsilon_{1} \in\left(0, \varepsilon_{1}^{*}\right)$ there exists $\varepsilon_{2}^{*}\left(\varepsilon_{1}\right)>0$ such that for each $\varepsilon_2 \in (0, \varepsilon_{2}^{*}\left(\varepsilon_{1}\right))$ $\ldots$ there exists $\varepsilon_{k}^{*}\left(\varepsilon_{k-1}\right)>0$ such that for each $\varepsilon_{k} \in$
$\left(0, \varepsilon_{k}^{*}\left(\varepsilon_{k-1}\right)\right)$ it holds:
    \begin{enumerate}
        \item (Semi-global stability) for each $R \geq \delta$, there exists $r > 0$, such that $\n{\phi(l, i)}_{\cA} \leq r \implies \n{\phi(t, j)}_{\cA} \leq R$ for $l + i \leq t + j$ and each solution $\phi$.
        \item (Practical attractivity)  for each $R, r$ that satisfy $\Delta \geq R \geq r \geq \delta$, there exists a period $T(r, R) \geq 0$, such that $\n{\phi(l, i)}_{\cA } \leq R \implies \n{\phi(t, j)}_{\cA} \leq r$ for all $t + j \geq T(r, R) + l + i$ and each solution $\phi$. \kraj
    \end{enumerate}
\end{defn}

\section{Singular perturbation theory for hybrid systems}
We consider two different system setups, with the first case featuring a hybrid reduced system and a continuous boundary layer system. In the second, both the reduced system and the boundary layer system are hybrid. Despite the different scenarios, we require similar assumptions in all configurations. Notably, we provide the most comprehensive coverage of the first case.
\subsection{Continuous boundary layer dynamics}
We consider the following hybrid dynamical system, denoted by $\mathcal{H}_1$:
\begin{subequations}\label{eq: singular main system}
\begin{align}
    \dot{x} &\in \m{I_{n_1} & 0 \\ 0 & \tfrac{1}{\varepsilon}I_{n_2}}F(x), &\text{if } x\in \cX_1 \times \cX_2, \label{eq: main system 1 flow}\\
    x^+ &\in G(x), &\text{if } x \in D_1 \times D_2, \label{eq: main system 1 jump}
\end{align} \label{eq: main system 1}
\end{subequations}
where $x \coloneqq \col{x_1, x_2} \in \mathcal{X}_1 \times \mathcal{X}_2 \subset \R^{n_1} \times \R^{n_2}$ are the system states, $\varepsilon > 0$ is small parameter used to speed up the $x_2$ dynamics, $\cX_1, D_1 \subset \mathcal{X}_1$, $\cX_2, D_2 \subset \mathcal{X}_2$ are flow and jump sets for the slow states $x_1$ and the fast states $x_2$, respectively. Other than $\varepsilon$, the system is implicitly parametrized by parameters $\beta, \gamma$ and $\tau$ i.e. $F = F_{\beta, \gamma, \tau}$ and $G = G_{\beta, \gamma, \tau}$. As it is common for hybrid dynamical systems, we postulate certain regularity assumptions that provide useful properties.
\begin{assum} \label{assum: hybrid conditions 1}
    The hybrid dynamical system in \nref{eq: main system 1} satisfies the basic regularity assumptions for hybrid systems \cite[Assum. 6.5]{goebel2012hybrid} for all parameters $\beta \in (0, \oln{\beta}], \gamma \in (0, \oln{\gamma}]$, $\tau \in (0, \oln{\tau}]$. The mapping $G$ satisfies item \cite[Assum. 6.5, A3]{goebel2012hybrid} also for $\beta = 0, \gamma = 0, \tau = 0$. Furthermore, all of systems's solutions are complete. \kraj
\end{assum}

Furthermore, we define two auxiliary systems in view of that in \nref{eq: singular main system}, the boundary layer system and the reduced system. The former, $\mathcal{H}_1^{\rho}$, for any given constant $\rho > 0$, is defined as
\begin{align}
    \dot{x} &\in \m{0 & 0 \\ 0 & I_{n_2}}F(x)\, & x\in \z{\z{\mathcal{A} + \rho \mathbb{B}}\cap \cX_1}\times \cX_2, \label{eq: bl system 1}
\end{align}
where $\mathcal{A} \subset \R^n$ is the equilibrium set of a reduced system, to be introduced later on. Furthermore, the system dynamics are parametrized by a small parameter $\beta$ which is used for tuning the desired convergence radius. In \nref{eq: bl system 1}, the dynamics of $x_1$ are frozen, i.e. $\dot{x}_1 = 0$, thus they approximate the behavior of those in \nref{eq: singular main system} when $\varepsilon > 0$ is chosen very small. Since the first state is constant, it is natural to assume that the equilibrium set, if it exists, contains all possible  $x_1$, i.e. the ones contained in the set $(\mathcal{A} + \rho \mathbb{B})\cap \cX_1$, and that for every $x_1$, there exists a specific set of equilibrium points $x_2$. We characterize this dependence with the ``steady-state" mapping $H$, and assume that it satisfies certain regularity properties \cite[Assum. 2]{sanfelice2011singular}, \cite[Assum. 2]{poveda2017framework}.

\begin{assum}\label{assum: singular main system basic conditions}
The set-valued mapping $H: \cX_1 \rightrightarrows \cX_2$, 
\begin{align}
    H(\oln{x}_1)\coloneqq \{ \oln{x}_2 \mid F(\oln{x}_1, \oln{x}_2) = 0\} \label{eq: assum steady state mapping}
\end{align}
is outer semicontinuous and locally bounded; for each $\oln{x}_1 \in \cX_1, H(\oln{x}_1)$ is a non-empty subset of  $\mathcal{X}_2$.  \kraj
\end{assum}
Now, we can define the complete equilibrium set of the system in \nref{eq: bl system 1}, the boundary layer manifold, as
\begin{align}
    \mathcal{M}_\rho \coloneqq \{ (x_1, x_2) \mid x_1 \in \z{\mathcal{A} + \rho \mathbb{B}} \cap \cX_1,\,  x_2 \in H(x_1)\}. \label{eq: mrho def}
\end{align} 
It is possible that the set $\mrho$ contains some unbounded states corresponding to the logic states or reference trajectories of the boundary layer system. We denote the bounded states with $x_2{'} \in \cX_2'$ and the unbounded states with $x_2^{\prime \prime} \in \cX_2^{\prime \prime}$, $\cX_2^{\prime} \times \cX_2^{\prime  \prime} = \cX_2$. Furthermore, we assume that these unbounded states only affect each other during jumps, and that the bounded states are a priori contained in a compact set. 

\begin{assum}\label{assum: jump mapping decomposition}
    The jump mapping $G$ in \nref{eq: main system 1 jump}, and the steady-state mapping $H$ in \ref{eq: assum steady state mapping} are decomposed as follows:
    \begin{align}
        G(x) &= \m{G_1(x_1, x_2') \\ G_2'(x)}, \label{eq: G' and G'' definition}\\
        H(x_1) &= H_1(x_1) \times \cX_2^{\prime \prime}, \label{eq: H definition}
    \end{align}
where $G_1: \cX_1 \times \cX_2' \rightrightarrows \cX_1 \times \cX_2'$, $G_2': \cX \rightrightarrows \cX_2^{\prime \prime}$, and $H_1: \cX_1  \rightrightarrows \cX_2'$.\kraj
\end{assum}
\begin{assum}\label{assum: boundness of x_2'}
    The set $\cX'_2$ in Assumption \ref{assum: jump mapping decomposition} is compact. \kraj
\end{assum}

Furthermore, we assume that the set $\mrho$ is SGPAS for boundary layer dynamics in Equation \nref{eq: bl system 1}.

\begin{assum} \label{assum: sp stability bl system}
The set $\mathcal{M}_\rho$ in \nref{eq: mrho def} is SGPAS as $\beta \rightarrow 0$ for the dynamics in \nref{eq: bl system 1}. Let $\Delta > \delta > 0$ be given by the definition of SGPAS. For every $\Delta > 0$, the corresponding Lyapunov function is given by
\begin{subequations}\label{eq: assum sp stability bl system 1}
\begin{align}
    &\alphadonjii \leq V_{2,\rho}(x) \leq \alphagornjii \\
    &\sup_{\begin{array}{l}
         \sm{f_1 \\ f_2}\in  F(x) 
    \end{array}}\vprod{\nabla V_{2,\rho}(x)}{\sm{0 \\ f_2}} \leq - \alpha_{2,\rho}(\n{x}_\mrho)\red
   &\text{for all }x\text{ such that }\n{x}_\mrho \geq \alpha_\beta(\beta), \\
   &\sup_{\begin{array}{l}
     \sm{f_1 \\ f_2}\in  F(x) 
    \end{array}}\vprod{\nabla V_{2,\rho}(x)}{\sm{0 \\ f_2}} \leq  \hat{\alpha}_{\beta}(\beta)\red
   &\text{for all }x\text{ such that }\n{x}_\mrho \leq \alpha_\beta(\beta) \\
   &\nabla V_{2,\rho}(x) = 0 \text{ for all }{x} \in \mrho,
\end{align}
\end{subequations}
where $\uad, \oad, \alpha_{2,\rho}, \alpha_\beta, \hat{\alpha}_{\beta}$ are functions of class $\mathcal{K}$, where $\alpha_{2,\rho}, \alpha_\beta$ are possibly parametrized by $\Delta$. Furthermore, for each compact set $K \in \cX_1$, there exists $M > 0$, such that
\begin{align}
\sup_{x \in K \times \cX_2} \n{V_{2, \rho}(x)} + \n{\nabla_{x_1}V_{2, \rho}(x)} \leq M. \kraj \label{eq: assum bl lyapunov V x2 bound}
\end{align}
\end{assum}

\begin{remark}
In Assumption \ref{eq: assum sp stability bl system 1}, we allow the set $\cX_2$ to be unbounded. Nevertheless, the Lyapunov function is assumed to take bounded values, as in \nref{eq: assum bl lyapunov V x2 bound}.\kraj
\end{remark}

On the other hand, since the $x_2$ dynamics are much faster than those of $x_1$ in \nref{eq: singular main system}, from the time scale of the latter, it seems that the $x_2$ dynamics are evolving on the manifold defined by the mapping $H$. To characterize this behaviour, we can define the reduced system $\mathcal{H}_1^{\text{r}}$ as:

\begin{subequations}
\begin{align}
    \dot{x}_1 &\in F_\text{r}(x_1)\,  & \text{ if } x_1 &\in \cX_1 \\
    x_1^+ &\in G_\text{r}(x_1)\,  & \text{ if }x_1 &\in D_1, \label{eq: reduced system 1 jumps}
\end{align} \label{eq: reduced system 1}
\end{subequations}
where $F_\text{r}(x_1) \coloneqq \overline{\operatorname{co}}\{ v_1 \mid (v_1, v_2) \in F(x_1, x_2), x_2 \in H(x_1) \}$, $G_\text{r}(x_1) \coloneqq \{ v_1 \mid (v_1, v_2) \in G(x_1, x_2), x_2 \in H(x_1)\}$. Furthermore, the system dynamics are parametrized by the parameter $\gamma$, which is used for the tuning of the convergence radius to the attractor set, and the parameter $\tau$ adjusts the minimum time interval between consecutive jumps, for \emph{almost all} jumps of the systems in \nref{eq: main system 1} (consequently also the reduced system in \nref{eq: reduced system 1}), as formalized next:
\begin{defn}[$\tau$-regular jump]
    A jump $j$ in a solution trajectory $\phi$ is a $\tau$-\emph{regular jump} if it occurs after an interval of flowing greater or equal than $\tau$, i.e. $\tau_j \coloneqq \sup \{\left|t-t^{\prime}\right|:(t, j - 1),\left(t^{\prime}, j - 1\right) \in \operatorname{dom} \phi\} \geq \tau$. Otherwise, the jump $j$ is called $\tau$-irregular. \kraj
\end{defn}
\begin{assum}\label{assum: regularity of jumps}
    Let $\phi$ be any solution of the system in \nref{eq: main system 1} with $\n{\phi(0, 0)}_{\cA\times\cX_2} \leq \Delta$. Then, there exists a finite number of jumps $N^*$ and finite time interval\, $T^*$, such that $\phi$ has at most $N^*$ $\underline{\sigma}(\tau)$-irregular jumps, and they all occur before $t \leq T^*$, where $\underline{\sigma}$ is a function of class $\mathcal{L}$, and $\tau$ is the parameter of the system.\kraj
\end{assum}

Differently from \cite{sanfelice2011singular}, where the reduced mapping is defined as $G_\text{r}(x_1) \coloneqq \{ v_1 \mid (v_1, v_2) \in G(x_1, x_2), x_2 \in \cX_2\}$, the mapping in \nref{eq: reduced system 1 jumps} only includes the jumps from the stead-state ``pairs" $(x_1, H(x_1))$ that belong to the manifold. Thus, our next assumption is weaker than \cite[Assum. 4]{sanfelice2011singular}, as it requires that the jumps stabilize the set $\mathcal{A}$ via a much more restricted set of dynamics. This is due to the fact that the reduced mapping $G_\text{r}$ does not contain all possible jumps from the set $D_1$, but only those from the boundary layer manifold $\mrho$. 
\begin{assum}\label{assum: reduced system stability 1}
The set $\mathcal{A}$ is SGPAS as $\gamma \rightarrow 0$ for the reduced system in \nref{eq: reduced system 1}. Let $\Delta > \delta > 0$ be given by the definition of SGPAS. For every $\Delta > 0$, the corresponding Lyapunov function is given by
\begin{subequations}\label{eq: assum reduced system stability 1}
\begin{align}
    &\alphadonji \leq V_1(x_1) \leq \alphagornji \\
    &\sup_{f_{1r} \in F_\text{r}(x_1)} \vprod{\nabla V_1(x_1)}{f_{1r}} \leq - \hat{\sigma}_{\tau}(\tau)\hat{\alpha}_{\gamma}(\gamma)\alpha_1\z{\n{x_1}_\cA} \label{eq: stdass reduced system lyapunov flows}\\
    &\sup_{g_{1r} \in G_\text{r}(x_1)} V_1(g_{r1}) - V_1(x_1) \leq - \hat{\alpha}_{\gamma}(\gamma)\alpha_1\z{\n{x_1}_\cA} \label{eq: stdass reduced sytem lyapunov jumps}\\
    &\text{for }\n{x_1}_{\mathcal{A}} \geq \alpha_\gamma(\gamma),
\end{align}
where $\uaj, \oaj, \alpha_1, \alpha_\gamma, \hat{\alpha}_{\gamma}$ are functions of class $\mathcal{K}$, where $\alpha_1, \alpha_\gamma$ are possibly parametrized by $\Delta$, and $\hat{\sigma}_{\tau}$ is a function of class $\mathcal{L}$.  
\end{subequations}
\kraj
\end{assum}

We claim that our original system in \nref{eq: main system 1} renders the set $\cA \times \cX_2$ practically attractive, if for almost all intervals of flow we allow the state of the system to converge to the neighborhood of the $\mrho$ manifold. The intuition is that in the neighborhood of the manifold, ``the jumps of the \emph{reduced system}" also contribute to the stabilization. 

\begin{theorem}\label{thm: basic case}
    Let Assumptions \ref{assum: hybrid conditions 1}---\ref{assum: reduced system stability 1} hold. Then the set $\mathcal{A} \times \cX_2$ is practically attractive as $(\gamma, \tfrac{1}{\tau}, \varepsilon, \beta) \rightarrow 0$ for the hybrid system in \nref{eq: main system 1}. \kraj
\end{theorem}

\begin{pf}
See Appendix \ref{proof: section theorem basic case}. \krajdokaz
\end{pf}

\begin{example}\label{example: attractivity}
Consider the hybrid dynamical system

\begin{subequations}\label{eq: basic example 1}

\begin{align}
    &\left\{ \begin{array}{rl}
         \dot{u} &= \gamma \max \{0, 1 - \tfrac{\n{u}}{R}\}  \\
         \dot{v} &= \tfrac{1}{\tau} \\
         \dot{x} &= -\tfrac{1}{\varepsilon}(x - u)
    \end{array}   \right. \red 
    &\text{if }(u, v, x) \in {[0, R] \times [0, 1] \times [0, R]};\\
    &\left\{ \begin{array}{rl}
         {u}^+ &= \tfrac{x}{2}  \\
         {v}^+ &= 0 \\
         {x}^+ &= R
    \end{array}   \right. \red 
    &\text{if }(u, v, x) \in {[0, R] \times \{1\} \times [0, R]},
\end{align}
    
\end{subequations}
where $\gamma, \tau, \varepsilon$ are tuning parameters, and $R > 0$ is the maximal trajectory radius. We show that the set $\{0\} \times [0, 1] \times [0, R]$ is practically attractive. First, we see that the boundary layer system reads as

\begin{align}
    &\left\{ \begin{array}{rl}
         \dot{u} &= 0 \\
         \dot{v} &= 0 \\
         \dot{x} &= -(x - u)
    \end{array}   \right. \red 
    &\text{if }(u, v, x) \in {[0, R] \times [0, 1] \times [0, R]}, \label{eq: example 1 bl system}
\end{align}
while the reduced system is given by 
\begin{subequations}
\begin{align}
    &\left\{ \begin{array}{rl}
         \dot{u} &= \gamma \max \{0, 1 - \tfrac{\n{u}}{R}\}  \\
         \dot{v} &= \tfrac{1}{\tau} \\
    \end{array}   \right. \red 
    &\text{if }(u, v) \in {[0, R] \times [0, 1]};\\
    &\left\{ \begin{array}{rl}
         {u}^+ &= \tfrac{u}{2}  \\
         {v}^+ &= 0 \\
    \end{array}   \right. \red 
    &\text{if }(u, v) \in {[0, R] \times \{1\}}.
\end{align}
\end{subequations}
Assumptions \ref{assum: hybrid conditions 1}---\ref{assum: regularity of jumps} are satisfied. Regarding Assumption \ref{assum: reduced system stability 1}, let the Lyapunov function of the reduced system be $V_1(u, v) = (2 - v) u^2$. It follows that 
\begin{align}
    &\dot{V}_1(u, v) \leq - \tfrac{1}{\tau}u^2 + 4\gamma a R, \red
    &V_1(u^+, v^+) - V_1(u, v) \leq - \tfrac{1}{2} u^2.
\end{align}

Since the reduced system satisfies Assumption \ref{assum: reduced system stability 1}, in view of Theorem \ref{thm: basic case}, practical attractivity is ensured. Unlike previous works \cite{wang2012averaging}, \cite{sanfelice2011singular}, and \cite{wang2012analysis}, our reduced jump mapping includes jumps only from the boundary layer, which allows us to establish stability results {using jumps}. In the aforementioned works, the reduced system jump mapping includes all possible jumps \cite[Equ. 13]{sanfelice2011singular}, \cite[Equ. 17]{wang2012analysis}, \cite[Equ. 13]{wang2012averaging}, and and for our example, it is given by $u^+ \in [-\tfrac{R}{2}, \tfrac{R}{2}]$. Thus, the assumption on the stability for reduced system dynamics \cite[Assum. 4]{sanfelice2011singular}, \cite[Thm. 2]{wang2012analysis}, \cite[Thm. 2]{wang2012averaging} does not hold. \kraj
\end{example}

We note that Theorem \ref{thm: basic case} gives us no guarantee on the stability of the state $x_1$, due to the fact that jumps can move the state arbitrarily far away from any set in $\cX_1$ (also seen in Example \ref{example: attractivity} for $u(0, 0) = 0, v(0, 0) = 0.99, x(0, 0) = R$). Under an additional assumption, it is possible to bound both the states $x_1$ and $x_2$ to a neighborhood of the set  $\mA \coloneqq \{ (x_1, x_2) \mid x_1 \in \mathcal{A}, x_2 \in H(x_1) \}$. 
\begin{assum}\label{assum: jumps in equilibrium}
The jump mapping $G$ in \nref{eq: main system 1 jump} is such that $G(\mA) \subset \mA$. \kraj
\end{assum}
Assumption \ref{assum: jumps in equilibrium} is sufficient to guarantee that for any neighborhood of the equilibrium set $\mA + \oln{r}\mathbb{B}$, there exists a neighborhood $\mA + \underline{r}\mathbb{B}$, such that jumps from the latter do not exit the former, i.e. $G(\mA + \underline{r}\mathbb{B}) \subset \mA + \oln{r}\mathbb{B}$. Lastly, we do not need to assume the compactness of the  set $\cX_2'$, as the distance from the set $\mA$ also bounds the values of the $x_2'$ state. 

\begin{theorem}\label{thm: constrainted jumps}
    Let Assumptions \ref{assum: hybrid conditions 1}---\ref{assum: jump mapping decomposition}, \ref{assum: sp stability bl system}---\ref{assum: jumps in equilibrium} hold. Then the set $\mA$ is SGPAS as $(\gamma, \tfrac{1}{\tau}, \varepsilon, \beta) \rightarrow 0$ for the hybrid system in \nref{eq: main system 1}. \kraj
\end{theorem}

\begin{pf}
See Appendix \ref{proof: section constrained jumps theorem}. \krajdokaz
\end{pf}

\begin{example}
We consider a hybrid dynamical system similar to one in \nref{eq: basic example 1}:
\begin{subequations}

\begin{align}
    &\left\{ \begin{array}{rl}
         \dot{u} &= \gamma\\
         \dot{v} &= \tfrac{1}{\tau} \\
         \dot{x} &= -\tfrac{1}{\varepsilon}(x - u)
    \end{array}   \right. \red 
    &\text{if }(u, v, x) \in {\R \times [0, 1] \times \R};\\
    &\left\{ \begin{array}{rl}
         {u}^+ &= \tfrac{x}{2}  \\
         {v}^+ &= 0 \\
         {x}^+ &= 2x
    \end{array}   \right. \red 
    &\text{if }(u, v, x) \in {\R \times \{1\} \times \R},
\end{align}
    
\end{subequations}
where $\gamma, \tau, \varepsilon$ are tuning parameters. Differently from \nref{eq: basic example 1}, the jump mapping is such that Assumption \ref{assum: jumps in equilibrium} is satisfied. Furthermore, as Theorem \ref{thm: constrainted jumps} does not require compactness of the set $\cX_2'$ in Assumption \ref{assum: boundness of x_2'}, the flow and jump sets are both unbounded. The boundary layer system has the same dynamics as the system in \nref{eq: example 1 bl system}, apart for the flow set which now reads as $\R \times [0, 1] \times \R$. The reduced system is given by 

\begin{subequations}
 \begin{align}\label{eq: example 2 dynamics}
&\begin{array}{l}
     \left\{ \begin{array}{rl}
         \dot{u} &= \gamma \\
         \dot{v} &= \tfrac{1}{\tau} \\
    \end{array}   \right.  \\
      \text{if }(u, v) \in {\R \times [0, 1]};
\end{array} \\
&\begin{array}{l}
     \left\{ \begin{array}{rl}
         {u}^+ &= \tfrac{u}{2}  \\
         {v}^+ &= 0 \\
    \end{array}   \right. \\
    \text{if }(u, v) \in {\R \times \{1\}}.
\end{array}
\end{align}   
\end{subequations}

Similarly to the previous example, all the Assumptions hold, thus due to Theorem \ref{thm: constrainted jumps},  the set $\{0\} \times [0, 1] \times \{0\}$ is SGPAS as $(\gamma, \tfrac{1}{\tau}, \varepsilon, \beta) \rightarrow 0$ for the dynamics in \nref{eq: example 2 dynamics}. Differently from \cite[Thm. 2]{wang2012averaging}, \cite[Thm. 1]{sanfelice2011singular}, and \cite[Thm. 2, Cor. 2 ]{wang2012analysis} where the \emph{fast states} are only a priori bounded to a compact set, here we can prove their convergence to the equilibrium set. \kraj
\end{example}

\subsection{Hybrid boundary layer dynamics}
Theorems \ref{thm: basic case} and \ref{thm: constrainted jumps} assume a lower limit on the time between \emph{all consecutive jumps} that occur in the system in \nref{eq: main system 1}. However, under certain conditions, it is possible to make a distinction between consecutive jumps of $x_1$, and the consecutive jumps of $x_2$. This is useful when the convergence of the boundary layer system is in fact driven by jumps in $x_2$, and imposing  a high lower limit on the period between consecutive jumps slows down convergence. Consider the following hybrid dynamical system, denoted with $\mathcal{H}_2$:
\begin{subequations}
\begin{align}
\dot{x} &\in \m{I_{n_1} & 0 \\ 0 & \tfrac{1}{\varepsilon}I_{n_2}} F(x), \text{ if } x\in \cX_1 \times \cX_2 \label{eq: main system 2 flow}\\
    x^+ &\in\left\{\begin{array}{l}
{\left[\begin{array}{l}
x_{1} \\
G_{2}(x)
\end{array}\right],\text{ if }  x \in \mathcal{X}_1 \times D_{2}} \\
{\left[\begin{array}{c}
G_{1}\left(x\right) \\
x_{2}
\end{array}\right], \text{ if }x \in D_{1} \times \mathcal{X}_{2}} \\
{\left[\begin{array}{l}
x_{1} \\
G_{2}(x)
\end{array}\right] \cup\left[\begin{array}{l}
G_{1}\left(x\right) \\
x_{2}
\end{array}\right], \text{ if }x \in D_{1} \times D_{2}}.
\end{array}\right.\label{eq: main system 2 jump}
\end{align}\label{eq: main system 2}
\end{subequations}
In this formulation, the distinction between the jumps of states $x_1$ and $x_2$ are highlighted, because during the jumps of $x_1$, $x_2$ stays constant, and vice versa. Furthermore, we define the boundary layer system, $\mathcal{H}_2^{\rho}$, as
\begin{subequations}
\begin{align}
    \dot{x} &\in \m{0 & 0 \\ 0 & I_{n_2}}F(x)\, & \text{if } x\in \cX_1 \times \cX_2, \\
    x^+ &\in \m{x_1 \\ G_2(x)}\, & \text{if } x \in \mathcal{X}_1 \times D_2,
\end{align} \label{eq: bl system 2}
\end{subequations}

and the reduced system, $\mathcal{H}_2^r$, as
\begin{subequations}
\begin{align}
    \dot{x}_1 &\in F_\text{r}(x_1)\,  & \text{ if }x_1 &\in \cX_1, \label{eq: restricted system 2 flows}\\
    x_1^+ &\in G_\text{r}(x_1)\,  &\text{ if }x_1 &\in D_1,
\end{align} \label{eq: reduced system 2}
\end{subequations}
where $F_\text{r}(x_1) \coloneqq \overline{\operatorname{co}}\{ v_1 \mid (v_1, v_2) \in F(x_1, x_2), x_2 \in H(x_1) \}$, $G_\text{r}(x_1) \coloneqq \{ v_1 \mid v_1 \in G_1(x_1, x_2),\, x_2 \in H(x_1)\}$. Differently from the boundary layer system in \nref{eq: bl system 1}, jumps are also included in this formulation, while the formulation of the reduced system is the same. Next, we pose analogous technical assumptions as for the system in \nref{eq: main system 1} and in turn provide results analogous to Theorem \ref{thm: constrainted jumps}.

\begin{assum} \label{assum: hybrid conditions 2}
    The hybrid dynamical system in \nref{eq: main system 2} satisfies the same conditions as in Assumption \ref{assum: hybrid conditions 1}.\kraj
\end{assum}
\begin{assum}\label{assum: G and H for system 2}
The jump mapping $G$ in \nref{eq: main system 2 jump}, and the steady-state mapping $H$ in \ref{eq: assum steady state mapping} are decomposed as in Equations \nref{eq: G' and G'' definition} and \nref{eq: H definition}. \kraj
\end{assum}

\begin{assum} \label{assum: sp stability bl system 2}
The set $\mathcal{M}_\rho$ is SGPAS as $\beta \rightarrow 0$ for the dynamics in \nref{eq: bl system 2}. Let $\Delta > \delta > 0$ be given by the definition of SGPAS. The corresponding Lyapunov function is given by \nref{eq: assum sp stability bl system 1}, with the additional equation
\begin{align*}
    \sup_{g_1 = x_1, g_2 \in G_2(x)}V_{2,\rho}(g) - V_{2,\rho}(x) \leq 0,
\end{align*}
and for each compact set $K \in \cX_1 \times \cX_2'$, there exists $M > 0$, such that
\begin{align}
\sup_{x \in K \times \cX_2^{\prime\prime}} \n{V_{2, \rho}(x)} + \n{\nabla_{x_1}V_{2, \rho}(x)} \leq M. \kraj \label{eq: assum bl lyapunov V x2' bound}
\end{align}
\end{assum}

\begin{assum}\label{assum: reduced system stability 2}
The set $\mathcal{A}$ is SGPAS as $\gamma \rightarrow 0$ for the reduced system in \nref{eq: reduced system 2}. Let $\Delta > \delta > 0$ be given by the definition of SGPAS. For every $\Delta > 0$, the corresponding Lyapunov function is given by Equation \nref{eq: assum reduced system stability 1} with the redefined mappings in \nref{eq: reduced system 2}.\kraj
\end{assum}

\begin{defn}[$\tau$-regular jump in $x_1$]
    A jump $j$ in a solution trajectory $\phi$ of the system in \nref{eq: main system 2} is a $\tau$-\emph{regular in $x_1$ jump} if it occurs after an interval of flowing in the $x_1$ state greater or equal than $\tau$, i.e. $\tau_j^1 \coloneqq \min \{\left|t-t^{\prime}\right| : \phi_1(t, j + 1) \in G_1(\phi(t, j)); \phi_2(t, j + 1) = \phi_2(t, j); (t, j), (t, j + 1) \in \operatorname{dom} \phi ; \phi_1(t', j' + 1) \in G_1(\phi(t', j')); \phi_2(t', j' + 1) = \phi_2(t', j'); (t', j'), (t', j' + 1) \in \operatorname{dom} \phi; j' > j \} \geq \tau$. Otherwise, if $\tau_j^1$ exists and $\tau_j^1 < \tau$,  jump $j$ is called $\tau$-irregular in $x_1$. \kraj
\end{defn}
\begin{assum}\label{assum: regularity of jumpsin x1}
    Let $\phi$ be any solution of the system in \nref{eq: main system 2} with $\n{\phi(0, 0)}_{\cA\times\cX_2} \leq \Delta$. Then, there exists a finite number of jumps $N^*$ and finite time interval\, $T^*$, such that $\phi$ has at most $N^*$ $\underline{\sigma}(\tau)$-irregular jumps in $x_1$, and they all occur before $t \leq T^*$, where $\underline{\sigma}$ is a function of class $\mathcal{L}$.\kraj
\end{assum}

\begin{corollary}\label{cor: hybrid bl system}
    Let Assumptions \ref{assum: singular main system basic conditions}, \ref{assum: jumps in equilibrium}---\ref{assum: regularity of jumpsin x1} hold. Then the set $\mA$ is SGPAS as $(\gamma, \tfrac{1}{\tau}, \varepsilon, \beta) \rightarrow 0$ for the hybrid system in \nref{eq: main system 2}. \kraj
\end{corollary}
\begin{pf}
The proof is analogous to the proofs of Theorems \ref{thm: basic case} and \ref{thm: constrainted jumps}. An equivalent for Lemma \ref{lemma: long enough jump times} can be constructed with jumps of the $x_2$ state. The rest of the proof is essentially the same.\krajdokaz
\end{pf}

\section{Illustrative example}

In \cite{stankovic2011distributed}, the issue of connectivity control was approached as a Nash equilibrium problem. In numerous practical situations, multi-agent systems are constructed with the goal of maintaining specific connectivity as a secondary objective in addition to their primary objective. In the subsequent discussion, we consider a comparable problem in which each agent is responsible for detecting an unknown signal source while also preserving a certain level of connectivity. Unlike \cite{stankovic2011distributed}, both the robots and the controllers have hybrid dynamics in our example.\\ \\
Consider a multi-agent system consisting of unicycle vehicles, indexed by $i \in \mathcal{I} \coloneqq\{1, \dots N\}$. Each agent is tasked with locating a source of a unique unknown signal. The strength of all signals abides by the inverse-square law, i.e. proportional to $1/r^2$. Therefore, the inverse of the signal strength can be used as a cost function. Additionally, the agents must not drift apart from each other too much, as they should provide quick assistance to each other in case of critical failure. This is enforced by incorporating the signal strength of the fellows agents in the cost functions. Thus, we design the cost functions as follows: 
\begin{align}
    \forall i \in \mathcal{I}: 
    h_i(u) = \|u_i - u_i^{\textup{s}} \|^2 + c\sum_{j \in \mathcal{I}_{-i}}\| u_i - u_j\|^2.
\end{align}
where $\mathcal{I}_{-i} \coloneqq \mathcal{I}\setminus \{i\}$, $c, b > 0$ and $u_i^{\textup{s}}$ represents the position of the source assigned to agent $i$. Goal of each agent is to minimize their cost function, and the solution to this problem is a Nash equilibrium. \\
\subsection{Unicycle dynamics}
As the unicycles are dynamical systems, a reference tracking controller is necessary in order to move them to the desired positions. In our example, let each agent implement a hybrid feedback controller similar to one in \cite{postoyan2015event} for trajectory tracking:
\begin{subequations}\label{eq: application unicycle dynamics}
\begin{align} 
&{\chi}_i^\textup{u} =  \col{{x}_i, {y}_i, \theta_i^e, {\tau}_i, \theta_i, \hat{v}_i, \hat{\omega}_i},\red
&\dot{\chi}_i^\textup{u}= F_i^\textup{u}(\chi_i^\textup{u}) \coloneqq\red
&\col{\hat{v}_i \cos \z{\theta_i}, \hat{v}_i \sin \z{\theta_i}, \omega_\text{r} - \hat{\omega}_i, \tfrac{1}{\sigma_i},  \hat{\omega}_i, 0, 0} \red \label{eq: new unicycle dynamics}
    &\text{ if } \chi_i^\textup{u} \in C_i^\textup{u} \coloneqq \R^3 \times [0, 1] \times \R^3,\\
&{\chi_i^\textup{u}}^+ = G_i^\textup{u}(\chi_i^\textup{u}) \coloneqq \col{x_i, y_i, \theta_i^e, 0, \theta_i, v_i, \omega_i} \red

    &\text{ if } \chi_i^\textup{u} \in D_i^\textup{u} \coloneqq \R^3 \times \{ 1 \} \times \R^3,
\end{align}
\end{subequations}
where $v_i = c_1(x_i^e - c_3 \omega_i y_i^e) - c_3 c_{2, i} (\omega_\textup{r} - \omega_i)y_i^e + c_3 \omega_i^2 x_i^e$, $x_i^e \coloneqq \cos(\theta_i)(u_i^1 - x_i) + \sin(\theta_i)(u_i^2 - y_i)$, $y_i^e \coloneqq -\sin(\theta_i)(u_i^1 - x_i) + \cos(\theta_i)(u_i^2 - y_i)$, $\theta_i^e = \theta_\textup{r} - \theta_i$, $\omega_i \coloneqq \omega_\textup{r} + c_{2, i} \theta_i^e$, $\dot{\theta}_\textup{r} = \omega_\textup{r} = const.$, $c_1, c_{2, i}, c_3> 0$ are tuning parameters, $\sigma_i$ is the sampling period parameter, $u_i^1$ and $u_i^2$ are the reference positions. Differently from \cite{postoyan2015event}, the jumps are triggered by a timer, and the reference trajectory is that of a unicycle with a fixed position $(u_i^1, u_i^2)$ and constant rotational velocity $\omega_\textup{r}$. Similarly to \cite[Lemma 4., Thm. 5]{postoyan2015event}, it is possible to prove that the dynamics in \nref{eq: application unicycle dynamics} render the set $\{\col{u_i^1, u_i^2, 0}\} \times \tilde{\mathcal{T}}_i \times \R^3$ SGPAS as $\sigma_i \rightarrow 0$.

\begin{theorem}\label{thm: unicycle refference convergance}
For $c_{2, i} = \sigma_i$, $c_3 = \tfrac{1}{3{{\omega}}_{\textup{r}}}$, $c_1 = \tfrac{1}{2c_3}$, the dynamics in \nref{eq: application unicycle dynamics} render the set $\{\col{u_i^1, u_i^2, 0}\} \times \tilde{\mathcal{T}}_i \times \R^3$ SGPAS as $\sigma_i \rightarrow 0$. \kraj
\end{theorem}

\begin{pf}
See Appendix \ref{proof: theorem of unicycle stampling}. \krajdokaz
\end{pf}

From the proof of Theorem \ref{thm: unicycle refference convergance}, it follows that system in  \nref{eq: application unicycle dynamics}, for all $i \in \mathcal{I}$, satisfies Assumptions \ref{assum: sp stability bl system 2}.\\

\subsection{Nash equilibrium seeking reference controller}
To steer the reference positions towards the Nash equilibrium, we implement the following asynchronous zeroth-order controller:

\begin{subequations} \label{eq: asyn algorithm}
\begin{align} 
    &\chi^\textup{c} = \col{\bfs{u}, \bfs{\xi}, \bfs{\mu}, \bfs{t}}, \red
    &\dot{\chi}^\textup{c} = F^\textup{c}(\chi^\textup{c}) \coloneqq \col{\bfs{0}, \bfs{0}, \bfs{0}, \bfs{\tau}^{-1}} \red
    &\text{ if } \chi^\textup{c} \in C^\textup{c} \coloneqq \R^{m} \times \mathcal{N} \times \mathbb{S}^{m} \times [0, 1]^N, \\
    &{\chi^\textup{c}}^+ = G^\textup{c}(\chi^\textup{c}), \text{i.e}.\red
    &\left\{ \begin{array}{ll}
        \bfs{u}^{+} = \bfs{u} - \alpha \beta S_x(\bfs{\tau})\bfs{\xi}\\ 
        \bfs{\xi}^{+} = \bfs{\xi} + \alpha S_x(\bfs{t})\z{2A^{-1} J(\bfs{x} + A\mathbb{D}\bfs{\mu})\mathbb{D}\bfs{\mu} - \bfs{\xi}} \\
        \bfs{\mu}^{+} = (I - S_{\mu}(\bfs{t}))\bfs{\mu} + S_{\mu}(\bfs{t}))\mathcal{R}\bfs{\mu} \\
        \bfs{t}^+ = (I - S_{\tau}(\bfs{t}))\bfs{t} 
             \end{array}\right. \red
    &\text{ if } \chi^\textup{c} \in D^\textup{c} \coloneqq \R^{m}  \times \mathcal{N} \times \mathbb{S}^{m}\times \mathcal{T}_\textup{R},
\end{align} 
\end{subequations}

where $\bfs{u} = \col{\z{u_i^1, u_i^2}_{i \in \mathcal{I}}}$ is used as the reference position for the systems in \nref{eq: application unicycle dynamics}$,  \bfs{\xi}$ is the collective filter state bound in a compact set $\mathcal{N} \subset \R^N$ chosen large enough to encompass all possible values of the state for all practical applications, $\bfs{\mu} \in \mathbb{S}^{2N}$ are oscillator states, $\bfs{t}$ are the timer states that control the sampling of each individual robot, $\bfs{\tau}^{-1} = \tau_0\col{\z{\tau_i^{-1}}_{i \in \mathcal{I}}}$ are the sampling periods that satisfy \cite[Assum. 9]{krilavsevic2023discrete}, $\bfs{x}$ are the positions of the unicycles, $\alpha, \beta > 0$ are small time-scale separation parameters, $\mathcal{R} \coloneqq \Diag{(\mathcal{R}_i)_{i \in \mathcal{I}}}$, $\mathcal{R}_i \coloneqq \Diag{\sm{\cos(\omega_i^j) & -\sin(\omega_i^j) \\ \sin(\omega_i^j) & \cos(\omega_i^j)}_{j \leq m_i}}$, $\omega_i^j > 0$ for all $i$ and $j$ are rotational frequencies and they satisfy \cite[Assum. 8]{krilavsevic2023discrete}, $\mathbb{D} \in \R^{2N \times 4N}$ is a matrix that selects every odd row from the vector of size $2N$, $a_i > 0$ are small perturbation amplitude parameters, $A \coloneqq \diag{(a_i)_{i \leq m}}$, $J(\bfs{x}) = \Diag{(J_i(x_i, \bfs{x}_{-i})I_{m_i})_{i \in \mathcal{I}}}$, $\mathcal{T} \subset [0, 1]^N$ is a closed invariant set in which all of the timers evolve and it excludes the initial conditions and their neighborhood for which we have concurrent sampling, $\mathcal{T}_\textup{R} \coloneqq \z{\cup_{i\in\mathcal{I}} [0, 1]^{i - 1} \times \{1\} \times [0, 1]^{N - i}}\cap\mathcal{T}$ is the set of timer intervals where one agent has triggered its sampling, $S_{x}: \mathcal{T} \rightarrow \R^{m \times m}$ and $S_{\tau}: \mathcal{T} \rightarrow \R^{N \times N}$ are continuous functions that output diagonal matrices with ones on the positions that correspond to states and timers of agents with $t_i = 1$, respectively, while other elements are equal to zero, when evaluating at $\bfs{t} \in \mathcal{T}_\textup{R}$.\\ \\

\subsection{The full system}

We define the collective state $\chi \coloneqq \col{\chi^\textup{c}, (\chi^\textup{u}_i)_{i \in \mathcal{I}}}$, collective flow map $F(\chi) \coloneqq \col{F^\textup{c}(\chi^\textup{c}), \tfrac{1}{\varepsilon}(F^\textup{u}_i(\chi^\textup{u}_i))_{i \in \mathcal{I}}}$, collective flow set $C\coloneqq {C^\textup{c} \times (C^\textup{u}_i)_{i \in \mathcal{I}}}$, collective jump map $G(\chi) \coloneqq \col{G^\textup{c}(\chi^\textup{c}), (G^\textup{u}_i(\chi^\textup{u}_i))_{i \in \mathcal{I}}}$, collective flow set $D\coloneqq (D^\textup{c} \times (C^\textup{u}_i)_{i \in \mathcal{I}}) \cup (C^\textup{c} \times (D^\textup{u}_i)_{i \in \mathcal{I}})$, and the equilibrium set $\cA_\chi \coloneqq \{\bfs{u}^*\} \times \mathcal{N} \times \mathbb{S}^N \times \mathcal{T} \times \{\col{\bfs{u}^*, \bfs{0}}\} \times [0, 1]^N \times \R^{3N}$. \\

\begin{figure}
    \centering
    \includegraphics{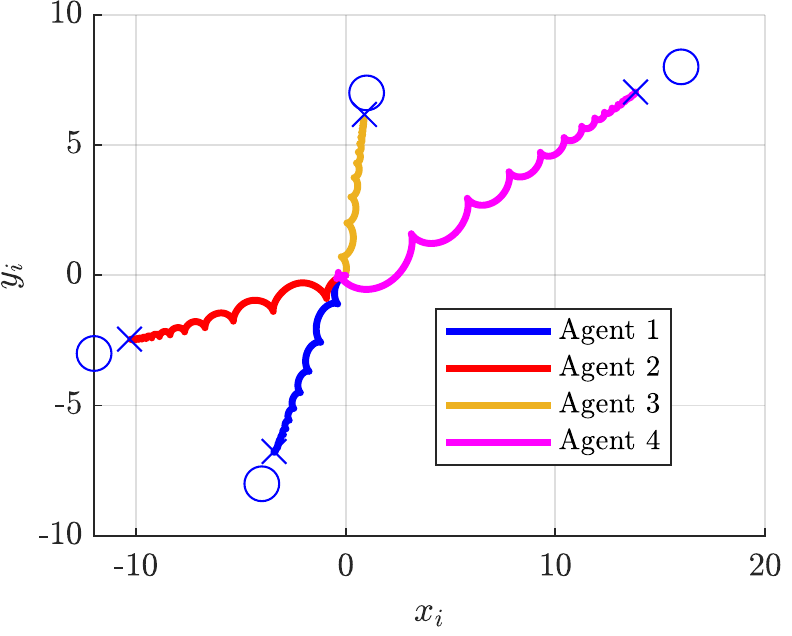}
    \caption{State trajectories in the $x-y$ plane. Circle symbols represent locations of the sources, while the $\times$ symbols represent locations of the NE.}
    \label{fig:application x phase plane}
\end{figure}

\begin{figure}
    \centering
    \includegraphics{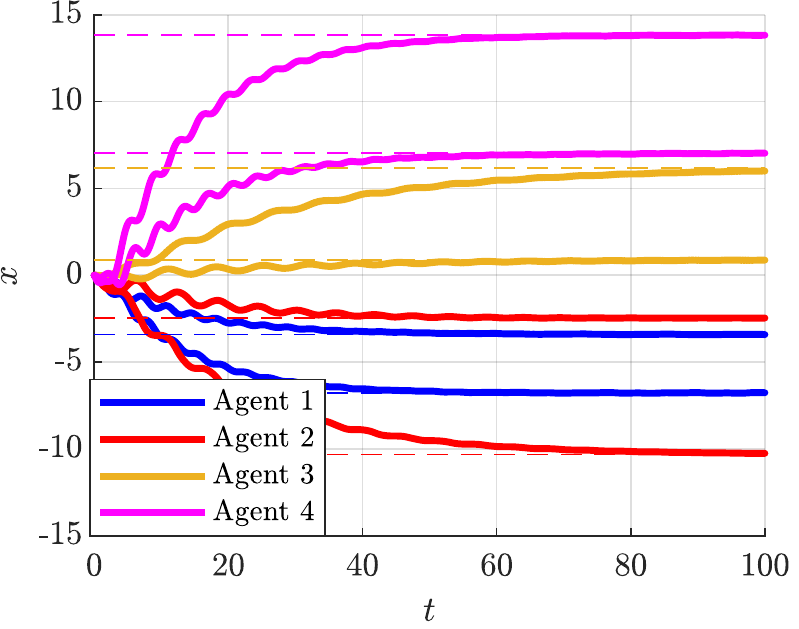}
    \caption{Time response of the unicycle position coordinates. The dashed lines correspond to the corresponding state of the Nash equilibrium.}
    \label{fig:application x time plot}
\end{figure}

We see that the steady state mapping is given by $H(\chi^\textup{c}) = \col{\bfs{u}, \bfs{0}} \times [0, 1]^N \times \R^{3N}$. Hence, the restricted system is equivalent to the one in \cite[Equ. 22]{krilavsevic2023discrete}. To show that Assumption \ref{assum: reduced system stability 2} is satisfied, we note that \cite[Thm. 1]{krilavsevic2023discrete} and \cite[Equ. E.10]{krilavsevic2023discrete} assure that the fully discrete-time zeroth-order variant of the algorithm in \cite[Equ. 22]{krilavsevic2023discrete}, has a Lyapunov function of the form 
\begin{align*}
    &\underline{\alpha}_{\textup{a}}\z{\n{\bfs{z} - \bfs{u}^*}} \leq V_{\textup{a}}(\bfs{z}) \leq  \oln{\alpha}_{\textup{a}}\z{\n{\bfs{z} - \bfs{u}^*}}  \\
    &V_{\textup{a}}(\bfs{z}^+) - V_{\textup{a}}(\bfs{z}) \leq - \hat{\alpha}_{\alpha}\z{{\alpha}}\alpha_{\textup{a}}\z{\n{\bfs{z} - \bfs{u}^*}} \\
    &\text{ for } \n{\bfs{z} - \bfs{u}^*} \geq \max\{\alpha_{\beta}(\beta), {\alpha}_{\alpha}(\alpha)\},
\end{align*}
where $\bfs{z} \coloneqq \bfs{u} + \bfs{\eta}$, and $\bfs{\eta}$ is a state of a bounded discrete system \cite[Equ. 7]{krilavsevic2023discrete}. For the sampled variant we have as our restricted system, we propose the following Lyapunov function 

\begin{align}
    V_1(\bfs{z}) \coloneqq & \tfrac{1}{2}\vprod{\bfs{1} - \bfs{t}}{\bfs{1}}\,\hat{\alpha}_{\alpha}\z{{\alpha}}\alpha_{\textup{a}} \z{\oln{\alpha}_{\textup{a}}^{-1}\z{\underline{\alpha}_{\textup{a}}\z{\n{\bfs{z} - \bfs{u}^*}}}} \red & + V_{\textup{a}}(\bfs{z}).
\end{align}
Hence, it holds 
\begin{align*}
    &\underline{\alpha}_{\textup{a}}\z{\n{\bfs{z} - \bfs{u}^*}} \leq V_1(\bfs{z}) \leq (\oln{\alpha}_{\textup{a}} + \alpha_{\textup{a}} \circ \oln{\alpha}_{\textup{a}}^{-1}\circ \underline{\alpha}_{\textup{a}} ) \z{\n{\bfs{z} - \bfs{u}^*}}  \\
    &\dot{V}_1(z) \leq - \tfrac{1}{2\tau_0} \sum_{i \in \mathcal{I}}\tau_i^{-1}\hat{\alpha}_{\alpha}\z{{\alpha}}\alpha_{\textup{a}} \z{\oln{\alpha}_{\textup{a}}^{-1}\z{\underline{\alpha}_{\textup{a}}\z{\n{\bfs{z} - \bfs{u}^*}}}}\\
    &V_1(\bfs{z}^+) - V_1(\bfs{z}) \leq - \tfrac{1}{2}\hat{\alpha}_{\alpha}\z{{\alpha}}\alpha_{\textup{a}}\z{\n{\bfs{z} - \bfs{u}^*}} \\
    &\text{ for } \n{\bfs{z} - \bfs{u}^*} \geq \max\{\alpha_{\beta}(\beta), {\alpha}_{\alpha}(\alpha)\},
\end{align*}
which satisfies Assumption \ref{assum: reduced system stability 2}. Furthermore, it is easy to show that Assumptions \ref{assum: singular main system basic conditions}, \ref{assum: jumps in equilibrium}, \ref{assum: hybrid conditions 2}, \ref{assum: G and H for system 2} hold as well. Since $\tau_0$ can be considered a tuning parameter for jump periods in the timers states $\bfs{t}$ in \nref{eq: asyn algorithm}, we can guarantee satisfaction of Assumption \ref{assum: regularity of jumpsin x1}. Hence, we satisfy all the Assumptions of the Corollary \ref{cor: hybrid bl system}, and for small enough parameters, the combined dynamics render the set $\cA_\chi$ SGPAS as $(\alpha, \beta, \max{\tau_i^{-1}}, \varepsilon, \max{\sigma_i})\rightarrow 0$. \\
For our numerical simulations, we choose the parameters: $u^s_1 = (-4, -8)$, $u^s_2 = (-12, -3)$, $u^s_3 = (1, 7)$, $u^s_4 = (16, 8)$, $(\sigma_1, \sigma_2, \sigma_3, \sigma_4) = \col{2, 3, 4, 2}\times 10^{-3}$, $c_1 = \tfrac{1}{3}$, $c_3 = 1.5$, $\alpha = 0.05$, $\beta = 0.003$, $c_{2, i} = \sigma_i$, $a_i = 0.1 $ for all $i$, $\bfs{t}(0, 0) = (0, 0.002, 0.004, 0.006)$, the perturbation frequencies $\omega_i^j$ were chosen as different natural numbers with added random numbers of maximal amplitude of 0.5, and the sampling of the Nash equilibrium seeking controller in \nref{eq: asyn algorithm} is five time slower than the sampling of the unicycle controller in \nref{eq: new unicycle dynamics}, i.e. $\bfs{\tau} = \col{1, 1.5, 2, 1} \times 10^{-2}$. \\
The numerical results are illustrated on Figures \ref{fig:application x phase plane} and \ref{fig:application x time plot}. We note that the trajectories converge to the neighborhood of the Nash equilibrium.


\section{Conclusion}
The application of singular perturbation theory can be extended to systems where the restricted system evolves on the boundary layer manifold through both flows and jumps. Moreover, by introducing some mild tehnical assumptions, one can show convergence of the fast state components towards a restricted attractor set that does not encompass the complete space of fast variables. With this theoretical extension, we can examine control systems that employ hybrid plants, along with controllers that are ``jump-driven" such as sampled controllers.
\endlinechar=13
\bibliographystyle{plain}        
\bibliography{biblioteka}  

\endlinechar=-1
\appendix

\section{Proof of Theorem \ref{thm: basic case}}\label{proof: section theorem basic case}
Let $\Delta > \delta > 0$ be given. We denote with $\rho$ the maximum distance to the equilibrium set $\cA \times \cX_2$ for trajectories starting in $(\cA + \Delta \mathbb{B})\times \cX_2$, which we characterize later on. Next, due to the fact that both $\cA \times \cX_2$ and $\mrho$ are unbounded in the dimensions corresponding to the same states, it follows that for any $\rho >0$, there exists a $P > 0$ such that $\n{x}_{\cA \times \cX_2} \leq \rho$ implies that $\n{x}_{\mrho} \leq P$. We consider the system in \nref{eq: singular main system} with restricted flow and jump sets:
\begin{subequations}\label{proof: restricted flow and jump sets}
\begin{align}
    C &\coloneqq \z{\z{\cA + \rho\mathbb{B}} \cap \cX_1} \times \cX_2 \\
    D &\coloneqq \z{\z{\cA + \rho\mathbb{B}} \cap D_1} \times D_2.
\end{align}
\end{subequations}

By plugging in $\Delta = P$ in Assumption \ref{assum: sp stability bl system}, and $\Delta = \rho$ in Assumptions \ref{eq: assum reduced system stability 1}, we construct the following Lyapunov function candidate:

\begin{align}
    V(x) = V_1(x_1) + \sqrt{\varepsilon} V_{2,\rho}(x). \label{proof1: lyapunov}
\end{align}

\subsection{Analysis of the jumps}
The Lyapunov function after jumps equals to 
\begin{align}
    V(g) = V_1(g_1) + \sqrt{\varepsilon} V_{2,\rho}(g),
\end{align}
where $\sm{g_1 \\ g_2} = g \in G(x)$. We prove the following Lemma:
\begin{lemma} \label{lemma: komsiluk u skokovima}
   Consider the hybrid system in \nref{eq: main system 1} with restricted flow and jump sets in \nref{proof: restricted flow and jump sets}, and let Assumptions \ref{assum: hybrid conditions 1}---\ref{assum: reduced system stability 1} hold. For every $e > 0$ and $\Delta > 0$, there exists $v > 0$, such that $\n{x}_{\mathcal{M}_\rho} \leq v$ implies that
    \begin{align*}
        &\sup_{\sm{g_1\\g_2}\in G(x)}V_1(g_1) -  V_1(x_1) \leq \red
        &\sup_{g_1^{\textup{r}} \in G_{\textup{r}}(x_1)}V_1(g_1^{\textup{r}}) -  V_1(x_1)+ \tfrac{e}{2}. \kraj
    \end{align*}
\end{lemma}

\begin{pf}
For the sake of contradiction, we assume that there exists $e > 0$ such it holds
\begin{align}
   &\sup_{\sm{g_1\\g_2}\in G(x)}V_1(g_1) -  V_1(x_1) > \red
   &\sup_{g_1^{\textup{r}} \in G_{\textup{r}}(x_1)}V_1(g_1^{\textup{r}}) -  V_1(x_1) + \tfrac{e}{2}, \red 
    &\forall x \in \z{\z{\mathcal{A} + \rho \mathbb{B}} \cap \cX_1}\times \cX_2 . \label{proof: lemma 1 counterstatement}
\end{align}
We define a sequence $\z{x^i}_{i \in \mathbb{N}} \in \mathcal{X}_1\times\mathcal{X}_2$ such that $\n{x^i}_\mrho \leq \tfrac{1}{i}$ and that it satisfies the inequality in Equation \nref{proof: lemma 1 counterstatement}. Let $x' \coloneqq \col{x_1, x_2'}$, and  $\mrho'$ be a projection of $\mrho$ onto the subspace of bounded states, $\cX_1 \times \cX_2'$. It holds that $\n{x^i}_\mrho \leq \tfrac{1}{i}$ implies $\n{{x^{\prime}}^i}_{\mrho{'}} \leq \tfrac{1}{i}$. Furthermore, it follows that the sequence ${x^{\prime}}^i$ is bounded. Due to Assumption \ref{assum: hybrid conditions 1}, we conclude that the sequence $\z{{g^{\prime}}^i}_{i \in \mathbb{N}}$, where  ${g^{\prime}}^i \in G_1({x^{\prime}}^i)$, is also bounded. Thus, due to the Weierstrass theorem, there exists a convergent subsequence that converges to the point $({x^\prime}^*, {g^\prime}^*)$, where ${x^\prime}^* \in \mrho'$. Next, due to the outer semi-continuity of the mappings $G$ and $H$, it holds that ${x_2^\prime}^* \in H_1(x_1^*)$, ${g^\prime}^* \in G_1({x^\prime}^*)$ and $g_1^* \in G_\text{r}(x_1^*)$. Therefore, it follows that 

\begin{align*}
    \sup_{\sm{g_1^*\\{g_2^\prime}^*}\in G_1(x_1^*, {x_2^\prime}^*)}&V_1(g_1^*) -  V_1(x_1^*) > \red
    \sup_{g_1^{\textup{r}} \in G_{\textup{r}}(x_1^*)}&V_1(g_1^{\textup{r}}) -  V_1(x_1^*)+ \tfrac{e}{2} \implies 0 > \tfrac{e}{2}\\
\end{align*}
which leads us to a contradiction and in turn proves the Lemma. \krajdokaz
\end{pf}

If $\varepsilon$ is chosen such that $\sqrt{\varepsilon} \leq \tfrac{e}{2 \overline{V}}$, where $\overline{V} \coloneqq \sup_{x \in \cX_1 \cap (\mathcal{A} + \rho \mathbb{B}) \times \cX_2, g \in G(x)}\n{V_{2,\rho}(g)}$,  then due to Lemma \ref{lemma: komsiluk u skokovima}, it holds that for any $e >0$ and $\Delta > 0$, there exist $v > 0$ and $\varepsilon^*$ such that for every $\varepsilon \in (0, \varepsilon^*)$, inequality  $\n{x}_\mrho \leq v$ implies that 
\begin{align}

        &\sup_{g \in G(x)} V(g) - V(x) \leq - \hat{\alpha}_{\gamma}(\gamma)\alpha_1(\n{x_1}_\cA)  + e\red
        &\text{for }\n{x_1}_{\mathcal{A}} \geq \alpha_\gamma(\gamma) \label{proof1: bounded jumping lyapunov}
\end{align}
The previous condition is always satisfied during jumps if $\n{x}_\mrho \leq v$ holds true before jumps. Thus, we have the following result:
\begin{lemma}\label{lemma: long enough jump times}
Consider the hybrid system in \nref{eq: main system 1} with restricted flow and jump sets in \nref{proof: restricted flow and jump sets}, and let Assumptions \ref{assum: hybrid conditions 1}, \ref{assum: singular main system basic conditions}, \ref{assum: sp stability bl system}  hold. Then, for every $v > 0$, and $\rho > 0$ there exists $\tau^* > 0$, $\varepsilon^*$ and $\beta^*$, such that for every $\varepsilon \in (0, \varepsilon^*)$, $\beta \in (0, \beta^*)$, it holds that if trajectory $x$ satisfies $\n{x(t, j)}_{\cA \times \cX_2} \leq \rho$ for all $t \in \dom(x(\cdot, j))$, then  $\n{x(t, j)}_\mrho \leq v$ for all  $t \in \dom(x(\cdot, j))$ such that $t \geq \tau^*$.\kraj
\end{lemma} 
\begin{pf}
From Assumption \ref{assum: sp stability bl system}, the derivative of the Lyapunov function candidate, for $\Delta = P$, reads as
\begin{align*}
     &\vprod{\nabla V_{2,\rho}(x)}{\sm{f_1 \\ \tfrac{1}{\varepsilon} f_2}} \leq -\tfrac{1}{\varepsilon}\alpha_{2,\rho}\z{\n{x}_\mrho} \\ &\quad \quad \quad \quad \quad \quad \quad \quad + \vprod{\nabla V_{2,\rho}(x)}{\sm{f_1 \\ 0}}\\
         &\text{for } \n{x}_\mrho \geq \alpha_\beta(\beta).\\
\end{align*}
We define the constant 
\begin{align*}
&\mu \coloneqq \sup_{x \in \cX_1 \cap \z{\mathcal{A} + \rho \mathbb{B}} \times \cX_2}\vprod{\nabla V_{2,\rho}(x)}{\sm{f_1 \\ 0}}.
\end{align*}
Then, the Lyapunov derivative is given by
\begin{align*}
    & \vprod{\nabla V_{2,\rho}(x)}{\sm{f_1 \\ \tfrac{1}{\varepsilon} f_2}} \leq -\tfrac{1}{\varepsilon}\alpha_{2,\rho}\z{\n{x}_\mrho} + \mu\\
    &\text{for } \n{x}_\mrho \geq \alpha_\beta(\beta).\\
\end{align*}

Let $\beta^* = \alpha_\beta^{-1}\z{\oad^{-1}\z{\uadz{v}}}$,\\ $m = \tfrac{1}{2}\alpha_{2,\rho}\z{\oad^{-1}\z{\uadz{v}}}$ and $\varepsilon^* = \tfrac{m}{\mu}$. It follows that for any time interval $(t, t + \tau)$ where only flowing occurred, for any $\varepsilon \in (0, \varepsilon^*)$, $\beta \in (0, \beta^*)$ it holds

\begin{align}
    &\dot{V}_{2,\rho}^P(x(t, j)) \leq - \mu \text{, for $\n{x}_\mrho \geq \oad^{-1}\z{\uadz{v}}$} \red
    &\int_{t}^{t + \tau} dV_{2,\rho}(x(t, j)) \leq - \mu \int_{t}^{t + \tau} dt \red
    &V_{2,\rho}(x(t + \tau, j)) \leq V_{2,\rho}(x(t, j)) - \mu\tau\red
    &\text{ for $\n{x}_\mrho \geq \oad^{-1}\z{\uadz{v}}$}. \label{eq: proof th1 bl layer lyapunov derivative}
\end{align}
As we assume $\n{x(t, j)}_\mrho \leq P$, from the bounds of the Lyapunov function in \ref{assum: sp stability bl system}, we have 

\begin{align}
    \n{x(t + \delta t, j)}_\mrho &\leq \uad^{-1}\z{\oadz{\n{x(t, j)}} - \mu\delta t}\red
    &\leq \uad^{-1}\z{\oadz{P} - \mu\delta t} \leq v. \label{eq: proof theorem sp attraction time bl system}
\end{align}
From the last inequality, it follows that $\tau \geq \tau^* \coloneqq \underline{\sigma}^{-1}
\z{\tfrac{\oadz{P} - \uadz{v}}{\mu}}$, which proves the Lemma. \krajdokaz

\end{pf}
It follows from \nref{proof1: bounded jumping lyapunov} and Lemmas \ref{lemma: komsiluk u skokovima} and \ref{lemma: long enough jump times} that for any $e > 0 $, $\Delta > 0$, there exist parameters $\varepsilon^*_1$, $\beta^*_1$, and $\tau^*$ such that for any $\varepsilon \in (0, \varepsilon^*_1)$, $\beta \in (0, \beta^*_1)$, if the time between consecutive jumps is larger than $\tau^*$, it holds that
\begin{align}
   &\sup_{g \in G(x)} V(g) - V(x)\leq  - \hat{\alpha}_{\gamma}(\gamma)\alpha_{1}\z{\n{x_1}_\cA} + e \label{eq: proof theorem sp jumps lyapunov}\\
   &\text{for }\n{x_1}_{\mathcal{A}} \geq \alpha_\gamma(\gamma).\nonumber
\end{align}
\subsection{Analysis of the flows}
The Lyapunov derivative is given by
\begin{align}
    &\sup_{\sm{f_1 \\ f_2} \in F(x)}\vprod{\nabla V(x)}{\sm{f_1 \\ \tfrac{1}{\varepsilon} f_2}}  = \red
    &\sup_{\sm{f_1 \\ f_2} \in F(x)}\z{\vprod{\nabla V_1(x)}{ f_1} + \sqrt{\varepsilon} \vprod{\nabla V_{2,\rho}(x)}{\sm{f_1 \\ \tfrac{1}{\varepsilon} f_2}}} \red
    & \leq \sup_{\sm{f_1 \\ f_2} \in F(x)} \left(\vprod{\nabla V_1(x)}{ f_1} + \left|\vprod{\nabla V_{2,\rho}(x)}{\sm{f_1 \\ 0}}\right| \right.\red 
    & \quad \left. +  \tfrac{1}{\sqrt \varepsilon}\vprod{\nabla V_{2,\rho}(x)}{\sm{0 \\ f_2}}\right) \red
    & \leq  -\hat{\sigma}_{\tau}(\tau)\hat{\alpha}_{\gamma}(\gamma)\alpha_1\z{\n{x_1}_\cA} + \mu(x) \red 
    & \quad + \sup_{\sm{f_1 \\ f_2} \in F(x)}\tfrac{1}{\sqrt \varepsilon}\vprod{\nabla V_{2,\rho}(x)}{\sm{0 \\ f_2}},\red
    &\text{ for }\n{x_1}_\cA \geq \alpha_\gamma(\gamma), \n{x}_\mrho \geq \alpha_\beta(\beta),
\end{align}
where $\mu(x) = -\sup_{f_1^\textup{r} \in F_\textup{r}(x_1)} \left(-\vprod{\nabla V_1(x)}{f_1^\textup{r}}\right)$ \\ $+ \sup_{\sm{f_1 \\ f_2} \in F(x), }\z{\vprod{\nabla V_1(x)}{ f_1} + \left|\vprod{\nabla V_{2,\rho}(x)}{\sm{f_1 \\ 0}}\right|}$.\\
Let $v > \alpha_{\beta}(\beta)$ be chosen arbitrarily and let 
\begin{align}
 \varepsilon^* = \tfrac{\underline{\alpha_{2\rho}}^2\z{v}}{\sup_{x \in \cX_1 \cap \z{\mathcal{A} + \rho \mathbb{B}} \times \cX_2} \n{\mu(x)}^2}. \label{proof1: varepsilon definition flows}
\end{align}
Then it holds that 
\begin{align}
    &\sup_{\sm{f_1 \\ f_2} \in F(x)}\vprod{\nabla V(x)}{\sm{f_1 \\ \tfrac{1}{\varepsilon} f_2}} \leq -\hat{\sigma}_{\tau}(\tau)\hat{\alpha}_{\gamma}(\gamma)\alpha_1\z{\n{x_1}_\cA}, \red 
    &\text{ for } \n{x}_\mrho \geq v, \n{x_1}_\cA \geq \alpha_\gamma(\gamma),  \label{proof1: lyapunov derivative outside v ball}\\
    &\sup_{\sm{f_1 \\ f_2} \in F(x)}\vprod{\nabla V(x)}{\sm{f_1 \\ \tfrac{1}{\varepsilon} f_2}} \leq -\hat{\sigma}_{\tau}(\tau)\hat{\alpha}_{\gamma}(\gamma)\alpha_1\z{\n{x_1}_\cA} \red 
    & + \mu(x) + \tfrac{1}{\sqrt{\varepsilon}}\hat{\alpha}_{\beta}(\beta), \text{ for } \n{x}_\mrho < v, \n{x_1}_\cA \geq \alpha_\gamma(\gamma),\label{proof1: lyapunov derivative inside v ball}
\end{align}
which is combined into 
\begin{align}
    &\sup_{\sm{f_1 \\ f_2} \in F(x)}\vprod{\nabla V(x)}{\sm{f_1 \\ \tfrac{1}{\varepsilon} f_2}} \leq -\hat{\sigma}_{\tau}(\tau)\hat{\alpha}_{\gamma}(\gamma)\alpha_1\z{\n{x_1}_\cA}\red
    & + \sup_{\n{x}_\mrho \leq v}\mu(x) + \tfrac{1}{\sqrt{\varepsilon}}\hat{\alpha}_{\beta}(\beta), \text{ for } \n{x_1}_\cA \geq \alpha_\gamma(\gamma).\label{proof1: lyapunov derivative}
\end{align}
The next Lemma shows that the positive terms in the Lyapunov derivative, with the proper choice of tuning parameters $\varepsilon$ and $\beta$, can be made arbitrarily small.
\begin{lemma}\label{lemma: komsiluk u tokovima}
Consider the hybrid system in \nref{eq: main system 1} with restricted flow and jump sets in \nref{proof: restricted flow and jump sets}, and let Assumptions \ref{assum: hybrid conditions 1}---\ref{assum: reduced system stability 1} hold. For every $e > 0$, $\Delta > 0$, there exists $\varepsilon^* > 0$, $\beta^*(\varepsilon) > 0$, such that for any $\varepsilon \in (0, \varepsilon^*)$, $\beta \in (0, \beta^*(\varepsilon))$, it holds that  $\sup_{\n{x}_\mrho \leq v}\mu(x) + \tfrac{1}{\sqrt{\varepsilon}}\hat{\alpha}_{\beta}(\beta) \leq e$. \kraj
\end{lemma}
\begin{pf}
We consider the following inequalities:
\begin{align}
    \sup_{\n{x}_\mrho \leq v}\mu(x) \leq \tfrac{e}{2} \label{proof1: lyapunov derivative remainder bound 1}\\
    \tfrac{1}{\sqrt{\varepsilon}}\hat{\alpha}_{\beta}(\beta) \leq \tfrac{e}{2}. \label{proof1: lyapunov derivative remainder bound 2}
\end{align}
If they hold, so does the inequality in the Lemma. The proof that inequality in \nref{proof1: lyapunov derivative remainder bound 1} can be made arbitrarily small by choice of small $v$, and hence smaller $\varepsilon^*$ due to \nref{proof1: varepsilon definition flows}, is analogous to the proof of Lemma \ref{lemma: komsiluk u skokovima}, and thus it is omitted. Then, to satisfy inequality \nref{proof1: lyapunov derivative remainder bound 2}, it is sufficient to have $\beta^* \coloneqq \min\z{\hat{\alpha}_{\beta}^{-1}\z{e \sqrt{\varepsilon}}, \alpha_\beta^{-1}(v)}$ and $\beta \in (0, \beta^*)$.\krajdokaz
\end{pf}

From Equations \nref{proof1: lyapunov derivative} and Lemma \ref{lemma: komsiluk u tokovima}, it follows that for any $e > 0$, $\Delta > 0$, there exists $\varepsilon^*_2$, $\beta^*_2(\varepsilon_2)$ such that for any $\varepsilon \in (0, \varepsilon^*_2)$ and $\beta_2 \in (0, \beta^*_2(\varepsilon_2))$, we have
\begin{align}
    &\sup_{\sm{f_1 \\ f_2} \in F(x)} \vprod{\nabla V(x)}{\sm{f_1 \\ \tfrac{1}{\varepsilon} f_2}} \leq -\hat{\sigma}_{\tau}(\tau)\hat{\alpha}_{\gamma}(\gamma)\alpha_1\z{\n{x_1}_\cA} \red 
    & + e, \text{ for }\n{x_1}_\cA \geq \alpha_\gamma(\gamma). \label{eq: proof theorem sp flows lyapunov}
\end{align}
\subsection{Complete Lyapunov analysis}

We denote by $\phi(t, j)$ a solution of the system that contains only $\underline{\sigma}(\tau)$ regular jumps. Let $\gamma$ be chosen so that $\gamma \in (0, \gamma^*)$, where $\gamma^* \coloneqq \min(\alpha_\gamma^{-1}(\oaj^{-1}(\tfrac{1}{2}\uajz{\delta})), \oln{\gamma})$. Next, $\eta$ is defined as $\eta \coloneqq \tfrac{1}{2}\min\{\uajz{\delta}, \alpha_1(\alpha_\gamma(\gamma^*)), 2\}$. Via Equation \nref{proof1: bounded jumping lyapunov} and Lemmas \ref{lemma: komsiluk u skokovima} and \ref{lemma: long enough jump times}, for $e = \hat{\alpha}_{\gamma}(\gamma)\eta$, we have $\tau^*$, $\varepsilon_1^*$, $\beta_1^*$. Next, we choose $\tau \in (0, \min(\tau^*, \oln{\tau}))$. From Equation \nref{proof1: lyapunov derivative} and Lemma \ref{lemma: komsiluk u tokovima} for $e = \hat{\sigma}(\tau)\hat{\alpha}_{\gamma}(\gamma)\eta$, we have $\varepsilon_2^*$, $\beta_2^*(\varepsilon)$. Finally, let 
\begin{align}
\varepsilon^*_3 = \tfrac{\eta^2}{{\sup_{x \in \cX_1 \cap (\mathcal{A} + \rho \mathbb{B}) \times \cX_2}}\alphagornjii^2}.\label{proof: varepsilon 3} 
\end{align}
We define $\varepsilon^* \coloneqq \min\{\varepsilon^*_1, \varepsilon^*_2\, \varepsilon^*_3, \oln{\varepsilon}\}$, $\beta^*(\varepsilon) \coloneqq \min\{\beta^*_1, \beta^*_2(\varepsilon), $ $\oln{\beta}\}$, and set the parameters as follows: $\varepsilon \in (0, \varepsilon^*)$, $\beta \in (0, \beta^*(\varepsilon))$. From Equations \nref{eq: proof theorem sp jumps lyapunov} and \nref{eq: proof theorem sp flows lyapunov}, it follows that

\begin{align}
    &V(\phi(t, j))+\sum_{i=0}^j \int_{t_i}^{t_{i+1}} \hat{\sigma}_{\tau}(\tau)\hat{\alpha}_{\gamma}(\gamma)\alpha_1\left(\n{\phi_1(s, i)}_{\mathcal{A}}\right) d s+ \red 
    &\sum_{i=1}^j \hat{\alpha}_{\gamma}(\gamma)\alpha_1\left(\n{\phi_1\left(t_i, i-1\right)}_{\mathcal{A}}\right) \leq V(\phi(0,0))\red
    &V(\phi(t, j))- V(\phi(0,0)) \leq \red & -(\hat{\sigma}_{\tau}(\tau)t + j)\hat{\alpha}_{\gamma}(\gamma)\z{\alpha_{1}\z{\n{x_1}_{\cA}} - \eta},\red
    &\text{ for }\n{\phi_1(t, j)}_\cA \geq \alpha_\gamma(\gamma).
\end{align}
As $\eta \leq \tfrac{1}{2}\alpha_1(\alpha_\gamma(\gamma))$ and $\alpha_\gamma(\gamma^*) < \alpha_\gamma(\gamma^*)$, we rewrite the last inequality as
\begin{align}
    &V(\phi(t, j))  \leq V(\phi(0,0)) - \tfrac{1}{2}(\hat{\sigma}_{\tau}(\tau)t + j)\alpha_1\z{\alpha_{\gamma}\z{\gamma^*}},\red
    &\text{ for }\n{\phi_1(t, j)}_\cA \geq \alpha_\gamma(\gamma^*), \label{eq: proof theorem sp lyapunov difference}
\end{align}
We can guarantee the decrease of the Lyapunov function up to the smallest Lyapunov level set that contains the set $(\cA + \alpha_\gamma(\gamma^*))\mathbb{B})\times \cX_2$. Via equation \ref{eq: proof theorem sp lyapunov difference}, we move onto proving semi-global boudness and practical attractivity of the equilibrium set. 
\subsubsection*{Semi-global boundness} 
 By definition, $\varepsilon \leq \varepsilon_3^*$. Thus, the upper and lower bound of the Lyapunov function candidate are given by 
\begin{align}
    \alphadonji \leq V(x) &\leq \alphagornji + \sqrt{\varepsilon} \alphagornjii \red 
    &\leqq \alphagornji + \eta\label{eq: proof theorem sp lyapunov bounds}
\end{align}

From \nref{eq: proof theorem sp lyapunov difference} and \nref{eq: proof theorem sp lyapunov bounds}, for $\underline{\sigma}(\tau)$-reggular trajectories $\phi$, it holds that 

\begin{align}
    \uajz{\n{\phi_1(t, j)}_\cA} \leq \oajz{\n{\phi_1(0, 0)}_\cA} + \eta, \red
    \text{ for }\n{\phi_1(t, j)}_\cA \geq \alpha_\gamma(\gamma^*). \label{eq: proof theorem sp semiglobal eq 1}
\end{align}

The maximal distance to the equilibrium of a $\underline{\sigma}(\tau)$-regular trajectory $\phi(t, j)$ starting at $(t, j) = (0, 0)$  is given in \nref{eq: proof theorem sp semiglobal eq 1} by $ \n{\phi_1(t, j)}_\cA \leq \rho_\textup{C}\z{\oajz{\n{\phi_1(0, 0)}_\cA}}$, where $\rho_\textup{C}: \R_+ \rightarrow \R_+, \rho_\textup{C}\z{a} \coloneqq \uaj ^{-1}\z{\oajz{a} + 1}$. By the outer semi-continuity and local boundness of the mapping $G$ in \nref{eq: G' and G'' definition} for all allowed sets of parameters, for  each $\underline{r} > 0$, there exists a $\oln{r}, r > 0$, such that $G_1\z{(\cA + \underline{r}\mathbb{B})\times\cX'_2} \subset G_1(\cA \times \cX'_2) + \oln{r}\mathbb{B} \subset (\cA + r\mathbb{B})\times \cX'_2$. Via this property, we define the mapping $\rho_\textup{D}: \R_+ \rightarrow \R_+, \rho_\textup{C}\z{\underline{r}} \coloneqq r$. Thus, for any initial condition such that $\n{\phi(0,0)}_{\cA \times \cX_2} \leq \Delta$, the maximal distance from the equilibrium set after $N$ irregular jumps, not necessarily consecutive jumps, is given by $\rho = \rho_\textup{C} \circ \rho_\textup{D}\circ \dots \circ \rho_\textup{C}(\Delta)$, where $\rho_\textup{D}$ repeats $N$ times, and $\rho_\textup{C}$ repeats $N + 1$ times. 

\subsubsection*{Semi-global stability for $\underline{\sigma}(\tau)$ trajectories} 

Let us consider trajectories after the $N$ irregular jumps. We show that for any $R \geq \delta$, there exists $r > 0$, such that $\n{\phi(l, i)}_{\cA \times \cX_2} \leq r \implies \n{\phi(t, j)}_{\cA \times \cX_2} \leq R$ for $l + i \leq t + j$. From \nref{eq: proof theorem sp semiglobal eq 1} and $\eta \leq \tfrac{1}{2}\uajz{\delta}$, it follows that 

\begin{align}
    &\n{\phi_1(l, i)}_\cA\leq \oaj^{-1}\z{\uajz{R} - \tfrac{1}{2}\uajz{\delta}}, \red
    &\text{ for }\n{\phi_1(t, j)}_\cA \geq \alpha_\gamma(\gamma^*). \label{eq: proof theorem sp semiglobal eq 2}
\end{align}

We note that the previous inequality holds up to the smallest radius of interest $r = \oaj^{-1}\z{\tfrac{1}{2}\uajz{\delta}}$, because $\gamma^* \leq \alpha_\gamma^{-1}(\oaj^{-1}\z{\tfrac{1}{2}\uajz{\delta}})$. Then, for any $R \geq \delta$, and all $\gamma \in (0, \gamma^*)$, we have $r(R) \coloneqq  \min\left\{\oaj^{-1}\z{\uajz{R} - \eta }, \Delta\right\}$. \\ \\

\subsubsection*{Practical attractivity}
Without of loss of generality, we assume $R, r$ are given so that $\Delta \geq R \geq r \geq \delta$. We have to show that there exists a period $T \geq 0$, such that $\n{\phi(0, 0)}_{\cA \times \cX_2} \leq R \implies \n{\phi(t, j)}_{\cA \times \cX_2} \leq r$ for all $(t, j)$ such that $t + j \geq T$.\\
Let $(l, i)$ be the hybrid time instant after the $N$ irregular jumps. By Assumption \ref{assum: regularity of jumps}, it holds $l + i \leq T^*$. Then, from \nref{eq: proof theorem sp lyapunov difference}, for $\phi(0, 0)$ replaced with $\phi(l, i)$, it follows

\begin{align}
    \uajz{r} \leq \oajz{\rho} - \tfrac{1}{2}\alpha_1\z{\alpha_{\gamma}\z{\gamma^*}}[(\hat{\sigma}_{\tau}(\tau)(t - l) + j - i)], 
\end{align}
Then, when we have
\begin{align}
    t + j \leq \tfrac{2(\oajz{\rho} - \uajz{r})}{\min\{\hat{\sigma}_{\tau}, 1\} \alpha_1\z{\alpha_{\gamma}\z{\gamma^*}}} + T^* = T(R, r).\label{eq: proof theorem sp attraction time}
\end{align}
\subsubsection*{Conclusion}
Our restricted system renders the set $\cA \times \cX_2$ practically attractive. Finally, to show the equivalence between the solutions of the original and restricted system, it is possible to use the same procedure as in \cite{sanfelice2011singular} after Equation (29). \krajdokaz

\section{Proof of Theorem \ref{thm: constrainted jumps}}\label{proof: section constrained jumps theorem}

Let $\Delta > \delta > 0$ be the parameters of semi-global practical stability. We denote with $\rho$ the maximum distance $\n{\phi(t, j)}_{\cA \times \cX}$ for trajectories starting in $\mA + \Delta \mathbb{B}$. On the other hand, as it is possible to a priori bound the distance $\n{\phi(t, j)}_{\mA}$ with $P > 0$, see Remark \ref{rem: bound P}, we ``redefine" the set $\cX_2'$ as a compact set $\z{\{x_2' \mid x_1 \in \cA, x_2' \in H_1(x_1)\} + P\mathbb{B}}\cap \cX'_2$. Now the same procedure as in proof of Theorem \ref{thm: basic case} can be repeated. From Equations \nref{eq: proof theorem sp attraction time bl system} and \nref{eq: proof theorem sp attraction time} we have
\begin{align}
 &\n{\phi(t, j)}_\mrho \leq v \text{ for all  } t \in \dom(\phi(\cdot, j)) \red
 &\text{ s.t. }t \geq \underline{\sigma}(\tau^*(P, v)), \label{proof2: bound x2}\\
 &\n{\phi(t, j)}_{\mathcal{A} \times \cX_2} \leq r  \text{ for all } (t, j) \in \dom(\phi) \red 
 &\text{ s.t. }t + j \geq T(\Delta, r),\label{proof2: KL bound x1}
\end{align}
 A key observation is that the intersection of sets $\mathcal{A} \times \cX_2$ and $\mrho$ gives us the set $\mA$, for which we want to prove stability. Distance to the set $\mathcal{A} \times \cX_2$ is given by \nref{proof2: KL bound x1}, and the distance to the set $\mrho$ is given by \nref{proof2: bound x2}. Thus, it is possible to quantify the distance to the equilibrium set using the distances of the latter two sets using the the following results:

\begin{lemma} \label{lemma: lower bound}
   Let $\mathcal{A}, \mathcal{B}$ be nonempty sets defined on a metric space, where at least one is bounded. Let their intersection $\mathcal{S}$ be nonempty. Then, for every $d > 0$, there exists $\underline{d} > 0$, such that $\n{x}_{\mathcal{A}} \leq \underline{d}$ and $\n{x}_{\mathcal{B}} \leq \underline{d}$ implies that $\n{x}_{\mathcal{S}} \leq d$. \kraj
\end{lemma}

\begin{pf}
Let us assume otherwise, i.e., there exists some $d > 0$ such that for any $\underline{d} > 0$ there exists $x$ such that it holds $\n{x}_{\mathcal{S}} > d$. Let us create a sequence of these points, $\z{x_i}_{i \in \mathbb{N}}$, such that $\n{x_i}_{\mathcal{A}} \leq \tfrac{1}{i}$ and $\n{x_i}_{\mathcal{B}} \leq \tfrac{1}{i}$. Because the sequence is bounded, there must exists a convergent subsequence. Let one such subsequence converge to $x^*$. Because of the continuity of the metric, it holds that $\n{x^*}_{\mathcal{A}} = 0$ and $\n{x^*}_{\mathcal{B}} = 0$. Thus, $x^* \in \operatorname{cl}\z{\mathcal{A}}$ and $x^* \in \operatorname{cl}\z{\mathcal{B}}$, or in other words $x^* \in \operatorname{cl}\z{\mathcal{A}} \cap \operatorname{cl}\z{\mathcal{B}} \equiv \operatorname{cl}\z{\mathcal{S}}$. Then it holds $\n{x^*}_{\mathcal{S}} = 0$, which is opposite of our assumption.\krajdokaz
\end{pf}

\begin{lemma}\label{lemma: upper bound}
   Let $\mathcal{A}, \mathcal{B}$ be nonempty sets defined on a metric space. Let their intersection $\mathcal{S}$ be nonempty and bounded. Then, for every $\overline{d} > 0$, there exists ${d} > 0$, such that $\n{x}_{\mathcal{S}} \leq d$ implies $\n{x}_{\mathcal{A}} \leq \overline{d}$ and $\n{x}_{\mathcal{B}} \leq \overline{d}$. \kraj
\end{lemma}

\begin{pf}
Let us assume otherwise, i.e there exists some $\overline{d} > 0$ such that for any $d > 0$, there exists $x$ such that $\n{x}_{\mathcal{A}} > \overline{d}$ and $\n{x}_{\mathcal{B}} > \overline{d}$. Let us create a sequence of these points,  $\z{x_i}_{i \in \mathbb{N}}$, such that $\n{x}_{\mathcal{A}} \leq \tfrac{1}{i}$. Because the sequence is bounded there must exist a convergent subsequence. Let one such subsequence converge to $x^*$. Because of the continuity of the metric, it holds that $\n{x^*}_{\mathcal{S}} = 0$. Thus, $x^* \in \operatorname{cl}\z{\mathcal{S}}$ or in other words $x^* \in \operatorname{cl}\z{\mathcal{S}} \equiv \operatorname{cl}\z{\mathcal{A}} \cap \operatorname{cl}\z{\mathcal{B}}$. Then it holds $\n{x^*}_{\mathcal{A}} = 0$ and $\n{x^*}_{\mathcal{B}} = 0$, which is opposite of our assumption.\krajdokaz
\end{pf}

\begin{remark}
Although we assume boundedness of some of the sets in Lemmas \ref{lemma: lower bound} and \ref{lemma: upper bound}, it is possible to prove the same results for the cases when $\cA$ and $\mathcal{B}$ are unbounded in the same dimensions, which is the case in our setup.
\end{remark}

\subsubsection*{Semi-global stability}

To prove practical stability, we show that for any $R \geq \delta$ there exists a neighborhood of the equilibrium, $\mA + r\mathbb{B}$, such that any trajectory initiated in neighborhood will stay inside the set $\mA + R\mathbb{B}$, for properly chosen parameters. But first, we prove a similar result for regular trajectories of the restricted system.\\

\begin{lemma}[Semi-global stability-like property]\label{lemma: stability for regular trajectories}
Consider the hybrid system in \nref{eq: main system 1} with restricted flow and jump sets in \nref{proof: restricted flow and jump sets}, and let Assumptions \ref{assum: hybrid conditions 1}, \ref{assum: singular main system basic conditions}, \ref{assum: jump mapping decomposition}, \ref{assum: sp stability bl system}, \ref{assum: reduced system stability 1} hold. Then, for every $\oln{v} \geq \delta$, there exists a $\underline{v} > 0$, and a set of tuning parameters $\varepsilon^*, \tau^*, \beta^*(\varepsilon), \gamma^*$, such that for all regular trajectories $\phi$ with $\n{\phi(0, 0)}_\mA \leq \underline{v}$ and $\varepsilon \in (0, \varepsilon^*), \tau \geq \tau^*, \beta \in (0, \beta^*(\varepsilon)), \gamma \in (0, \gamma^*)$, it holds that $\n{\phi(t, j)}_\mA \leq \oln{v}$ for all $(t, j) \in \dom\z{\phi}$. \kraj
\end{lemma}
\begin{pf}

\emph{Sketch of the proof}\\
First, we find $\hat{v} > 0$ such that any trajectory initiated in $\mA + \hat{v}\mathbb{B}$, stays in $\mA + \overline{v}\mathbb{B}$ during flows. Then we find $\tilde{v} > 0$ such that jumps from $\mA + \tilde{v}\mathbb{B}$ will end in $\mA + \hat{v}\mathbb{B}$. Next, we find $\underline{v}$ such that any trajectory initiated in $\mA + \underline{v}\mathbb{B}$, stays in $\mA + \tilde{v}\mathbb{B}$ during flows. Finally, we choose $\varepsilon$, $\tfrac{1}{\tau}$, $\beta$, $\gamma$ small enough such that all trajectories end up in $\mA + \underline{v}\mathbb{B}$ before jumps. \\ \\

Consider the following system of implications:

\begin{align*}
   &\quad\,\,\,\n{\phi(0, 0)}_{\mA} {\leq} {\underline{v}} \overset{(1)}{\Rightarrow} \begin{array}{l}
         \n{\phi(0, 0)}_{\cA \times \cX_2} \leq \underline{u} \\
         \n{\phi(0, 0)}_{\mrho} \leq \underline{u}
    \end{array} \\
    &\overset{(2)}{\Rightarrow} \begin{array}{l}
         \n{\phi(t, 0)}_{\cA \times \cX_2} \leq \tilde{u} \\
         \n{\phi(t, 0)}_{\mrho} \leq \tilde{u}
    \end{array} \overset{(3)}{\Rightarrow} \n{\phi(t, 0)}_{\mA} \leq \tilde{v} \\
    &\overset{(4)}{\Rightarrow} \n{\phi(t, 1)}_{\mA} \leq \hat{v} \overset{(5)}{\Rightarrow} \begin{array}{l}
         \n{\phi(t, 1)}_{\cA \times \cX_2} \leq \hat{u} \\
         \n{\phi(t, 1)}_{\mrho} \leq \hat{u}
    \end{array} \\
   &\overset{(6)}{\Rightarrow} \begin{array}{l}
         \n{\phi(l, 1)}_{\cA \times \cX_2} \leq \overline{u} \\
         \n{\phi(l, 1)}_{\mrho} \leq \overline{u}
    \end{array} \overset{(7)}{\Rightarrow} \n{\phi(l, 1)}_{\mA} \leq \overline{v}. \label{proof: semi global stability of regular trajectories}
\end{align*}
Implication (7) follows from Lemma \ref{lemma: lower bound}, while implication (6) follows from Equations \nref{eq: proof theorem sp semiglobal eq 2} with $\hat{u}^1 \leq \oaj^{-1}\z{\uajz{\tfrac{1}{2}\overline{u}}}$, Equation \nref{eq: proof th1 bl layer lyapunov derivative} and $\hat{u}^1 = \oad^{-1}\z{\uadz{\overline{u}}}$, and $\hat{u} = \min\{\hat{u}^1, \hat{u}^2\}$; Implication (5) follows from Lemma \ref{lemma: upper bound}; Implication (4) proceeds from \cite[Lemma 5.15]{goebel2012hybrid}, outer semicontinuity, local boundedness of the mapping $G$, Assumption \ref{assum: jumps in equilibrium}, thus for every $\hat{v} > 0$, there exists a $\tilde{v} \leq \hat{v}$ such that $G(\mA + \tilde{v}\mathbb{B}) \subset \mA + \hat{v}\mathbb{B}$; Implication (3) follows from Lemma \ref{lemma: lower bound}, while implication (2) follows from Equations \nref{eq: proof theorem sp semiglobal eq 2} with $\underline{u}^1 = \oaj^{-1}\z{\uajz{\tfrac{1}{2}\tilde{u}}}$, Equation \nref{eq: proof th1 bl layer lyapunov derivative} and $\underline{u}^2 = \oad^{-1}\z{\uadz{\tilde{u}}}$ and $\underline{u} = \min\{\underline{u}^1, \underline{u}^2\}$; Implication (1) follows from Lemma \ref{lemma: upper bound}. To satisfy the inequalities in Equations \nref{eq: proof theorem sp semiglobal eq 2}, \nref{eq: proof th1 bl layer lyapunov derivative}, let $\eta = \tfrac{1}{2}\min\{\uajz{\underline{u}}, \alpha_1\z{\oaj^{-1}\z{\tfrac{1}{2}\uajz{\underline{u}}}}, 2\}$,  $\gamma^* \coloneqq $\\ $\alpha_\gamma^{-1}\z{\oaj^{-1}\z{\tfrac{1}{2}\uajz{\underline{u}}}}$, $\gamma \in (0, \gamma^*)$; Via Equation \nref{proof1: bounded jumping lyapunov} and Lemmas \ref{lemma: komsiluk u skokovima} and \ref{lemma: long enough jump times}, for $e = \hat{\alpha}_{\gamma}(\gamma)\eta$, we have $\tau^*$, $\varepsilon_1^*$, $\beta_1^*$. Next, we choose $\tau \in (0, \min(\tau^*, \oln{\tau}))$. From Equation \nref{proof1: lyapunov derivative} and Lemma \ref{lemma: komsiluk u tokovima} for $e = \hat{\sigma}(\tau)\hat{\alpha}_{\gamma}(\gamma)\eta$, we have $\varepsilon_2^*$, $\beta_2^*(\varepsilon)$. Finally, let $\varepsilon_3^*$ be defined as in \nref{proof: varepsilon 3}. We define $\varepsilon^* \coloneqq \min\{\varepsilon^*_1, \varepsilon^*_2\, \varepsilon^*_3, \oln{\varepsilon}\}$, $\beta^*(\varepsilon) \coloneqq \min\{\beta^*_1, \beta^*_2(\varepsilon), $ $\oln{\beta}\}$, and set the parameters as follows: $\varepsilon \in (0, \varepsilon^*)$, $\beta \in (0, \beta^*(\varepsilon))$.  \\
Furthermore, as Equation \nref{eq: proof theorem sp semiglobal eq 2} holds for jumps and flows, it follows that $\n{\phi(l, 1)}_{\cA \times \cX_2} \leq \tilde{u}$, and due to Lemma \ref{lemma: long enough jump times}, it holds that $\n{\phi(l, 1)}_{\mrho} \leq \tilde{u}$ for $l \geq \min_t \dom\z{\phi(\cdot, 1)} + \underline{\sigma}(\tau^*)$. Thus, it is possible to follow the same reasoning with implications (3) to (7) for the next, and all the following  regular jumps, which proves our Lemma.  \krajdokaz
\end{pf}
\begin{remark}\label{rem: bound P}
We can ``reverse" Lemma \ref{lemma: stability for regular trajectories} so that we claim that for every $\underline{v} > 0$, there exists a $\oln{v} > 0$ that satisfies the same inequality. Then, by doing an inverse procedure of the proof of stability, we can derive $P > 0$ such that $\n{\phi(0, 0)}_\mA \leq \Delta$ implies that $\n{\phi(t, j)}_\mA \leq P$. These bounds depend on the proprieties of mapping $G$ and the lower and upper bounds of the Lyapunov functions, thus can be computed a priori. \kraj
\end{remark}
Let $N$ be the number of irregular jumps for the given $\Delta$. Via Lemma \ref{lemma: stability for regular trajectories}, for $\oln{v} = R$, we find $\oln{r}_N = \underline{v}$ and parameters $\varepsilon^*_N, \tau^*_N, \beta^*_N(\varepsilon), \gamma^*_N$. Then, from \cite[Lemma 5.15]{goebel2012hybrid}, outer semicontinuity, local boundedness of the mapping $G$, Assumption \ref{assum: jumps in equilibrium}, we can find $\underline{r}_N$ such that $G(\mA + \underline{r}_N\mathbb{B}) \subset \mA + \oln{r}_N\mathbb{B}$. Then again we use Lemma \ref{lemma: stability for regular trajectories}, with $\oln{v} = \underline{r}_N$, to find $\underline{v} = \oln{r}_{N-1}$ and parameters $\varepsilon^*_{N - 1}, \tau^*_{N - 1}, \beta^*_{N - 1}(\varepsilon), \gamma^*_{N - 1}$. These steps are repeated until we reach the first jump. Then, we use Lemma \ref{lemma: stability for regular trajectories}, for $\oln{v} = \underline{r}_1$ to find $r = \underline{v}$ and parameters $\varepsilon^*_{0}, \tau^*_{0}, \beta^*_{0}(\varepsilon), \gamma^*_{0}$. We note that for $\varepsilon^* \coloneqq \min\{ \varepsilon^*_0, \dots, \varepsilon^*_N \}$, $\tau^* \coloneqq \max\{ \tau^*_0, \dots, \tau^*_N \}$, $\beta^*(\varepsilon) \coloneqq \min\{ \beta^*_0(\varepsilon), \dots, \beta^*_N(\varepsilon) \}$ and $\gamma^* \coloneqq \min\{ \gamma^*_0, \dots, \gamma^*_N \}$, all the inequalities hold.

\noindent\emph{\textbf{Practical attractivity}} \\ \\

Without of loss of generality, we assume $R, r$ are given so that $\Delta \geq R \geq r \geq \delta$. Let $(l, i)$ be a hybrid time instant after the $N$ irregular jumps. By Assumption \ref{assum: regularity of jumps}, it holds $l + i \leq T^*$. Furthermore, Lemma \ref{lemma: stability for regular trajectories} gives us $\underline{r} = \underline{v}$ for $\oln{v} = \delta$, and the corresponding tuning parameters $\varepsilon^*, \tau^*, \beta^*(\varepsilon), \gamma^*$. From the definition of parameters in Lemma \ref{lemma: stability for regular trajectories}, it follows that the Lyapunov derivatives and differences for functions in Equation \nref{proof1: lyapunov} and Assumption \ref{assum: sp stability bl system}, are defined for $\n{\phi(t, j)}_{\cA \times \cX_2} \geq \underline{u}, \n{\phi(t, j)}_{\mrho} \geq \underline{u}$, where $\underline{u}$ is given in the the system of implications in \nref{proof: semi global stability of regular trajectories}. Thus Equations \nref{proof2: bound x2} and \nref{proof2: KL bound x1} guarantee that the trajectories eventually enter and stay in $\tilde{u}$ neighborhoods before jumps, for $v = r = \underline{u}$. And from our practical-stability result, it follows that the trajectory stays in the $\oln{r}$ neighborhood \krajdokaz

\section{Proof of Theorem \ref{thm: unicycle refference convergance}}\label{proof: section of unicycle theorem}\label{proof: theorem of unicycle stampling}

Similarly to \cite[Equ. 13]{postoyan2015event}, let the Lyapunov function candidate be given by 
\begin{align}
    &V_i(q_i)\coloneqq\tfrac{1}{2} \z{{x}_i^e - c_3{\omega}_i y_i^e}^2 +\tfrac{1}{2} {y^e_i}^2+\tfrac{1}{2} {\theta^e_i}^2,
\end{align}

where $q_i = \col{{x}_i, {y}_i, \theta_i^e, {\tau}_i, \theta_i, \hat{v}_i, \hat{\omega}_i}$. First, we characterize the upper and lower bounds of the Lyapunov function candidate. It holds
\begin{align*}
    &V_i(q_i) = \tfrac{1}{2}{{x}_i^e}^2 -c_3{\omega}_i x_i^e y_i^e + \tfrac{1}{2}c_3^2{\omega}_i^2{y_i^e}^2 +\tfrac{1}{2} {y^e_i}^2+\tfrac{1}{2} {\theta^e_i}^2 \\
    &\geq \tfrac{1}{2}\, {{x}_i^e}^2 \z{1 - \gamma} + \tfrac{1}{2}c_3^2{\omega}_i^2 {y_i^e}^2 \z{1 - \tfrac{1}{\gamma}} + \tfrac{1}{2} {y^e_i}^2+\tfrac{1}{2} {\theta^e_i}^2\\
    &{\geq} \tfrac{1}{4}{{x}_i^e}^2 + \tfrac{1}{2}\z{1 - c_3^2{\omega}_i^2}{y^e_i}^2 + \tfrac{1}{2} {\theta^e_i}^2 \\
    &\geq \tfrac{1}{4}\z{x_i - u_i^1}^2 + \tfrac{1}{4}\z{y_i - u_i^2}^2+\tfrac{1}{4} {\theta^e_i}^2 = \tfrac{1}{4}\n{r_i}^2,
\end{align*}
where the second line follows from $ab \leq \tfrac{\gamma}{2}a^2 + \tfrac{1}{2\gamma}b^2$, third line follows from $\gamma = \tfrac{1}{2}$, in the forth line we assume that $c_3 \leq \tfrac{\sqrt{2}}{2{\omega}_i}$, and $r_i \coloneqq \col{x_i - u_i^1, y_i - u_i^2, \theta_i - \theta_{\textup{r}}}$. Furthermore, for the upper bound we have
\begin{align*}
    V_i(q_i) &\leq {{x}_i^e}^2 + c_3^2{\omega}^2{y_i^e}^2 +\tfrac{1}{2} {y^e_i}^2+\tfrac{1}{2} {\theta^e_i}^2\\
    &\leq {{x}_i^e}^2 + {y^e_i}^2 + {\theta^e_i}^2 = \n{r_i}^2,
\end{align*}
where the second line follows from the former assumption on constant $c_3$. Thus, the bound on the Lyapunov function are given by.
\begin{align*}
    \tfrac{1}{4}\n{r_i}^2 \leq V_i(q_i) \leq \n{r_i}^2.
\end{align*}
The Lyapunov derivative is bounded similarly to \cite[Equ. 14]{postoyan2015event}:
\begin{align*}
    \vprod{\nabla V_i(q_i)}{f_i(q_i)} = -\Sigma_i(q_i) + \Lambda_i(q_i),
\end{align*}
where $\Sigma_i(q_i) \coloneqq c_1\z{{x}^e_i - c_3{\omega}_i y^e_i}^2 + c_2  {\theta^e_i}^2 + c_3 \omega_i^2 {y^e_i}^2$ and $\Lambda_i(q_i)\coloneqq\z{{x}^e_i - c_3{\omega}_i y^e_i}\left(e^w_i y^e_i-e^v_i-c_3 y^e_i e^w_i c_2\right.$ $\left.+c_3 w e^\omega_i x^e_i\right)-y^e_i e^w_i x^e_i-\theta^e_i e^w_i$, with $e^w_i:=$ $\hat{\omega}_i-\omega_i$, $e^v_i = \hat{v}_i-v_i$, and $v_i \coloneqq c_1(x_i^e - c_3 \omega_i y_i^e) - c_3 c_2 (\omega_\textup{r} - \omega_i)y_i^e + c_3 \omega_i^2 x_i^e$. To characterize the convergence rate, we upper bound $\Sigma(q_i)_i$ as follows:

\begin{align*}
    &\Sigma_i(q_i) = c_1\z{{x}^e_i - c_3{\omega}_i y^e_i}^2 + c_2  {\theta^e_i}^2 + c_3 \omega_i^2 {y^e_i}^2\\
    &\geq c_1{{x}^e_i}^2\z{1 - \gamma} + c_1 c_3^2\omega_i^2 {y^e_i}^2\z{1 - \tfrac{1}{\gamma}} + c_3\omega_i^2 {y^e_i}^2 + c_2  {\theta^e_i}^2 \\
    &\geq \tfrac{1}{2}c_1{{x}^e_i}^2 + \z{c_3\omega_i^2 - c_1 c_3^2\omega_i^2}{y^e_i}^2 + c_2 {\theta^e_i}^2\\
    &\geq \tfrac{1}{2}c_1{{x}^e_i}^2 + \tfrac{1}{2}c_3\omega_i^2 {y^e_i}^2 +  c_2 {\theta^e_i}^2 \\
    &\geq \underline{c}\n{r_i}^2,
\end{align*}
where the third line follows for $\gamma = \tfrac{1}{2}$, in fourth line we assume $c_1 \leq \tfrac{1}{2 c_3}$, $\underline{c}\coloneqq \min\{\tfrac{1}{2}c_1, c_2, \tfrac{1}{2}c_3\omega_i^2\}$  Now, we write the Lyapunov derivative as 
\begin{align*}
    \vprod{\nabla V_i(q_i)}{f_i(q_i)} \leq - \underline{c}V_i(q_i) + \Lambda_i(q_i).
\end{align*}
As the jumps restart $\Lambda_i(q_i)$ to $0$, the jumps of the Lyapunov given by
\begin{align*}
    V(q_i^+) - V(q_i) \leq 0.
\end{align*}
Let $\Delta > \delta > 0$ be the parameters of the semi-global practical stability. If our Lyapunov derivative is negative on the desired domain, it follows that
\begin{align*}
    \tfrac{1}{4}\n{r_i(t, j)}^2 \leq V_i(q_i(t, j))\leq V_i(q_i(0, 0)) \leq \n{r_i(0,0)}^2,
\end{align*}
Thus for any initial condition with $\n{r_i(0,0)} \leq \Delta$, it holds that $\n{r_i(t, j)} \leq 2\Delta$. Using the previous bound, we can estimate the minimal and maximal value of $\omega_i$ as 
\begin{align}
    \underline{\omega}_i &\coloneqq \min \omega_i = \omega_{\textup{r}} - 2 c_2 \Delta \\
    \overline{\omega}_i &\coloneqq \max \omega_i = \omega_{\textup{r}} + 2 c_2 \Delta.  
\end{align}
Hence, we choose $c_2 = \tfrac{\omega_{\textup{r}}}{4\Delta}$ to ensure $\underline{\omega}_i$ is positive, $c_3 = \tfrac{1}{3{{\omega}}_{\textup{r}}} \leq \tfrac{\sqrt{2}}{3{{\omega}}_{\textup{r}}}$, and $c_1 = \tfrac{1}{2c_3}$. \\
As $\Lambda_i(q_i)$ is differentiable and all of its variables and their derivatives are bounded, we can approximate it with a constant $M$ and write the derivative as 

\begin{align*}
    \vprod{\nabla V_i(q_i)}{f_i(q_i)} &\leq - \underline{c}V_i(q_i) + M \tilde{\tau}_i \\
    &\leq - \underline{c}V_i(q_i) + M \sigma_i 
\end{align*}
Parameter $\sigma_i$ can be made arbitrarily small, thus enabling arbitrarily close convergence to the equilibrium point. As there is a constant time between jumps, semi-global practical stability follows for $c_1  = \sigma_i \rightarrow 0$ \cite[Cor. 8.7]{sanfelice2011singular}.
\end{document}